\documentclass{ieeetj}
\usepackage{cite}
\usepackage{amsmath,amssymb,amsfonts}
\usepackage{algorithmic}
\usepackage{algorithm}
\usepackage{graphicx,color}
\usepackage{textcomp}
\usepackage{hyperref}
\usepackage{lmodern}
\hypersetup{}
\usepackage{algorithm,algorithmic}
\def\BibTeX{{\rm B\kern-.05em{\sc i\kern-.025em b}\kern-.08em
    T\kern-.1667em\lower.7ex\hbox{E}\kern-.125emX}}
\AtBeginDocument{\definecolor{tmlcncolor}{cmyk}{0.93,0.59,0.15,0.02}\definecolor{NavyBlue}{RGB}{0,86,125}}
\def\OJlogo{\vspace{-4pt}}
\def\seclogo{\vspace{10pt}}
\def\authorrefmark#1{\ensuremath{^{\textbf{#1}}}}
\usepackage{comment}

\newcommand{\bsm}{\mathbf{m}}

\newcommand{\bsy}{\mathbf{y}}
\newcommand{\bsw}{\mathbf{w}}
\newcommand{\bsz}{\mathbf{z}}
\newcommand{\bstheta}{\boldsymbol{\theta}}

\onecolumn

\usepackage[table]{xcolor}

\definecolor{my_gray}{HTML}{DBDBDB}

\usepackage{booktabs, makecell, multirow, array}
\usepackage{makecell}

\newcolumntype{a}{>{\columncolor{my_gray}}c}

\usepackage{amsthm}

\begin{document}
\receiveddate{}
\reviseddate{}
\accepteddate{}
\publisheddate{}
\currentdate{}
\doiinfo{}

\markboth{Benchmarking Diffusion Annealing-Based Bayesian Inverse Problem Solvers}{Scope Crafts and Villa}

\title{ Benchmarking Diffusion Annealing-Based Bayesian Inverse Problem Solvers}

\author{Evan Scope Crafts\authorrefmark{1}, Member, IEEE, and Umberto Villa\authorrefmark{1,2}, Member, IEEE}
\affil{Oden Institute for Computational Engineering and Sciences, The University of Texas at Austin, Austin, TX USA 78712}
\affil{Dept of Biomedical Engineering, The University of Texas at Austin, Austin, TX USA 78712}
\corresp{Corresponding author: Umberto Villa (email: uvilla@austin.utexas.edu).}
\authornote{This work was supported by the National Institute of Biomedical Imaging and Bioengineering of the National Institutes of Health under award numbers R01EB031585 and R01EB034261. Code is available at: \url{https://doi.org/10.5281/zenodo.14908136}. Datasets containing the results of the numerical studies can be found at: \url{ https://doi.org/10.7910/DVN/0L5KGB}. }

\begin{abstract}

In recent years, the ascendance of diffusion modeling as a state-of-the-art generative modeling approach has spurred significant interest in their use as priors in Bayesian inverse problems. However, it is unclear how to optimally integrate a diffusion model trained on the prior distribution with a given likelihood function to obtain posterior samples. While algorithms developed for this purpose can produce high-quality, diverse point estimates of the unknown parameters of interest, they are often tested on problems where the prior distribution is analytically unknown, making it difficult to assess their performance in providing rigorous uncertainty quantification. Motivated by this challenge, this work introduces three benchmark problems for evaluating the performance of diffusion model based samplers. The benchmark problems, which are inspired by problems in image inpainting, x-ray tomography, and phase retrieval, have a posterior density that is analytically known. In this setting, approximate ground-truth posterior samples can be obtained, enabling principled evaluation of the performance of posterior sampling algorithms. This work also introduces a general framework for diffusion model based posterior sampling, Bayesian Inverse Problem Solvers through Diffusion Annealing (BIPSDA). This framework unifies several recently proposed diffusion-model-based posterior sampling algorithms and contains novel algorithms that can be realized through flexible combinations of
design choices. We tested the performance of a set of BIPSDA algorithms, including previously proposed
state-of-the-art approaches, on the proposed benchmark problems. The results provide insight into the strengths and limitations of existing diffusion-model based posterior samplers, while the benchmark problems provide a testing ground for future algorithmic developments.
\end{abstract}


\maketitle

\section{INTRODUCTION}

Inverse problems, which aim to estimate an unknown parameter of interest from observed data (measurements), provide a principled way to integrate data with existing scientific knowledge  \cite{GhattasWillcox21,mueller2012linear}. However, many important inverse problems throughout the sciences are \textit{ill-posed}  \cite{engl1996regularization,tarantola2005inverse}---the measurements alone do not contain enough information to uniquely and stably estimate the parameter of interest \cite{tenorio2017introduction}. In the Bayesian inverse problem framework, this difficulty is addressed through the integration of prior knowledge regarding the parameter of interest, which enables both point estimation of the unknown parameter and rigorous uncertainty quantification \cite{kaipio2006statistical, stuart2010inverse}. This can facilitate risk-informed analysis and decision making, which is particularly important in the context of safety-critical systems (e.g., medical imaging \cite{prince2006medical} or earthquake detection \cite{slingerland2011mathematical} systems), where accurate risk assessments can save lives. 

The Bayesian framework requires knowledge of the prior probability distribution of the unknown parameter. However, for many problems of interest (e.g., image reconstruction) the distribution of the unknown parameters has a complex structure that is not easily captured by hand-crafted prior distributions (e.g., Gaussian, total variation, or sparsity-inducing priors). Using an inaccurate prior can lead to biased point estimates of the unknown parameter and incorrect predictions of the parameter uncertainty. 

The success of generative modeling \cite{chan2024tutorial,ruthotto2021introduction} in capturing the structure of complex high-dimensional probability distributions in recent years has spurred considerable interest in the use of generative modeling to overcome this limitation \cite{bora2017compressed}. In particular, over the last four years the rise to prominence of diffusion modeling \cite{song2019generative, ho2020denoising, song2020score} as a state-of-the-art generative modeling approach has been coupled with significant interest in their use as prior distributions in Bayesian inference. Here a common strategy is to obtain posterior samples by pretraining a diffusion model on samples from the prior distribution and integrating the pretrained model with the likelihood function (which is assumed to be known analytically) at inference time \cite{jalal2021robust, song2022solving, chung2022improving, chung2023diffusion, wu2024principled, darassurvey}. This general strategy has the advantage of enabling the same diffusion model to be used in conjunction with many different data acquisition designs without retraining, and has shown considerable promise. However, there is no consensus on how to optimally integrate the diffusion model and the likelihood function. While several algorithms have been proposed for this purpose, test-bed problems with no analytically-known ground-truth prior are often used to illustrate the performance of such methods. In this setting, the performance of the algorithms is often assessed using sample quality metrics like peak signal-to-noise ratio, as well as heuristic arguments regarding the sample diversity. This makes it difficult to assess the accuracy of these methods in capturing the global structure of the posterior. 

In this paper, we provide a general framework, Bayesian Inverse Problem Solvers through Diffusion Annealing (BIPSDA), which---by abstracting the algorithmic components of the Diffusion Annealing Posterior Sampling (DAPS) method described in \cite{zhang2024improving}---provides an unified formulation for developing and analyzing annealing-based solvers for Bayesian inverse problems with priors implicitly defined by a diffusion model. 
The framework generalizes and extends two recently proposed algorithms, DAPS \cite{zhang2024improving} and DiffPIR  \cite{zhu2023denoising}, for solving Bayesian inverse problems with diffusion models. These two algorithms have achieved strong performance on a number of canonical imaging reconstruction problems, including non-linear problems such as phase retrieval that were previously considered too difficult for diffusion-model-based solvers \cite{zhang2024improving}. In the unified framework, the DAPS and DiffPIR algorithms can be recovered through specific design choices, while new algorithms can be unveiled through exploration of the rich algorithmic design space. An original contribution of our work is the use of randomize-then-optimize (RTO) techniques, originally proposed in \cite{bardsley2014randomize,Kainan2018randomized,ba2022randomized} for approximate sampling from posterior distributions, to solve a sampling subproblem that arises in the BIPSDA framework. Here the RTO technique transforms this subproblem into an optimization problem (as in DiffPIR) that can be solved with ``off-the-shelf'' deterministic numerical optimization methods, while still accurately accounting for measurement noise statistics (as in the DAPS method).

We systematically evaluated the performance of algorithms in the proposed framework, including the DAPS and DiffPIR algorithms and the novel RTO-based algorithms, on a set of model problems. Each model problem uses a Gaussian mixture prior. This choice of prior has two key advantages. First, under this choice the posterior density is known and approximate ground-truth posterior samples can be obtained, enabling rigorous analysis of the performance of the BIPSDA algorithms. Second, under this choice the noisy prior scores, a key component of diffusion models that is learned from data, can be formed analytically. This allows us to decompose the error in the posterior sampler into two components: error inherent to the algorithm, and error due to incorrect modeling of the prior distribution. For the likelihood functions, we considered benchmark problems inspired by classic image restoration/reconstruction problems: simple linear inpainting problems (in both low and high noise regimes), as well as nonlinear x-ray tomography and phase retrieval based problems. 

In each experiment, the performance of the BIPSDA algorithms was evaluated by comparing the samples produced by each algorithm with the approximate ground-truth samples using four different metrics: the central moment discrepancy (CMD) \cite{zellinger2017central}, maximum mean discrepancy (MMD) \cite{gretton2012kernel}, and the errors in both the predicted posterior mean and pointwise variance. This enables principled evaluation of the ability of various BIPSDA algorithms to capture both the local and global posterior structure. 

The remainder of this paper is organized as follows. In Section \ref{sec:background}, we provide relevant background on Bayesian inverse problems and diffusion models. In Section \ref{sec:proposed} we introduce the proposed BIPSDA framework, discuss its relationship to previously proposed algorithms, and provide examples of algorithms that can be realized within the framework. Section \ref{sec:numerical_staples} provides details regarding the numerical studies, while results are shown in Section \ref{sec:results}. A discussion of the results and the conclusion is given in Section \ref{sec:discussion_conc}. 

\section{BACKGROUND}
\label{sec:background}

In this section, we provide relevant background on Bayesian inverse problems and diffusion models. We also provide a brief overview of a popular class of algorithms, referred to here as hijacking algorithms, for solving Bayesian inverse problems with diffusion models. This overview serves to both motivate the proposed framework and to introduce algorithmic ideas relevant to the present work.

\subsection{BAYESIAN INVERSE PROBLEMS}

In inverse problems, the relationship between the measurements $\bsy \in \mathbb{R}^K$ and the unknown variable of interest $\bsm \in \mathbb{R}^D$ is captured by the likelihood function $\pi_{\mathrm{like}}(\bsy \mid \bsm)$, where here and throughout the remainder of this work we assume the measurements and the unknown variable lie in finite-dimensional Euclidean spaces. Concretely, under the assumption of additive Gaussian noise, the measurements can be written as 
\begin{equation}
\label{eq:additive-gaussian}
    \bsy = f(\bsm) + \bsz, \quad \bsz \sim \mathcal{N}(\bsz; \mathbf{0}, \boldsymbol{\Sigma}_{\bsz}),
\end{equation}
where $f: \mathbb{R}^D \to \mathbb{R}^K$ is known as the forward model, and $\boldsymbol{\Sigma}_{\bsz} \in \mathbb{R}^{K \times K}$ is the covariance matrix of the noise distribution. The corresponding likelihood function is given as $\pi_{\mathrm{like}}(\bsy \mid \bsm) = \mathcal{N}(\bsy; f(\bsm), \boldsymbol{\Sigma}_{\bsz})$.

By Bayes' Theorem, the posterior distribution $\pi_{\mathrm{post}}(\bsm \mid \bsy)$ of the unknown variable is related to the likelihood function and prior density function $\pi_{\mathrm{pr}}(\bsm)$ by the following expression:
\begin{equation}
   \pi_{\mathrm{post}}(\bsm \mid \bsy) \propto \pi_{\mathrm{like}}(\bsy \mid \bsm) \; \pi_{\mathrm{pr}}(\bsm). 
\end{equation}
The goal of Bayesian inverse problems is to characterize the posterior distribution. In particular, depending on the application, various quantities, such as the mean, covariance, or higher-order moments of the posterior, may be of interest. While for certain choices of priors and likelihood functions (e.g., linear-Gaussian likelihood and Gaussian prior), these quantities can be computed in closed form, in general this is not tractable, and in practice Monte Carlo methods are often used to estimate the quantities of interest from samples. 

\subsection{DIFFUSION MODELS}

Broadly, the goal of generative modeling can be described as follows: given a set of samples $\mathcal{M} = \{ \bsm_i \}_{i=1}^{N_s}$ from a distribution of interest with density function $\pi_{0}(\bsm)$, obtain a set of samples from $\pi_{0}(\bsm)$ that are not in $\mathcal{M}$. In diffusion modeling, these samples are obtained by reversing a predefined noising process \cite{ho2020denoising, song2020score}. This noising process can be written as a stochastic differential equation (SDE). In particular, the variance exploding (VE) form of the diffusion SDE can be written as follows \cite{song2020score}:
\begin{equation}
\label{eq:vanilla_sde}
d \bsm(t) = \sqrt{2 \dot{\sigma}(t)\sigma(t)} \; d\bsw(t), \quad t \in [0, T], 
\end{equation}
where $\sigma(t)$ is a predefined noise schedule, $\bsw(t)$ is the Wiener process, and $\bsm(0) \sim \pi_{0}(\bsm)$. Under this SDE, the conditional distribution of $\bsm(t)$ given $\bsm(0)$ is given by the \textit{noising distribution} $\pi_{t \mid 0}(\bsm(t) \mid \bsm(0)) = \mathcal{N}(\bsm(t); \bsm(0), \sigma^2(t) \mathbf{I})$, where here and throughout the remainder of this work we have assumed $\sigma(0) = 0$ for simplicity. This implies that solving the forward SDE is equivalent to sampling from a Gaussian. Further, for large enough $T$, the distribution $\pi_{T}(\cdot)$ will be approximately Gaussian.

To sample from $\pi_{0}(\bsm)$, we can first sample from the Gaussian approximation of $\pi_{T}(\bsm(T))$ and then reverse the noising process in \eqref{eq:vanilla_sde}. Here there are two main methods to reverse the noising process \cite{song2020score}: an SDE-based approach and an ordinary differential equation (ODE) based approach. The SDE-based approach leverages the remarkable fact \cite{anderson1982reverse} that \eqref{eq:vanilla_sde} admits a time-reversal that matches the marginal distributions $\pi_t (\bsm(t))$. This SDE takes the following form:
\begin{align}
d\bsm(t) &= -2 \dot{\sigma}(t)\sigma(t) \nabla_{\bsm(t)} \log \pi_{t}(\bsm(t)) \; dt  \nonumber \\
& \quad \quad + \sqrt{2 \dot{\sigma}(t)\sigma(t)} \; d\bsw(t). \label{eq:reverse_SDE}
\end{align}
The ODE based approach is based on the fact that there exists an ODE that has the same marginal distributions $\pi_{t}(\bsm(t))$ as the forward and reverse SDEs. This ODE, known as the probability flow ODE, has the following form:
\begin{align}
d\bsm(t) &= - \dot{\sigma}(t)\sigma(t) \nabla_{\bsm(t)} \log \pi_{t}(\bsm(t)) \; dt. 
\label{eq:prob_flow_ODE}
\end{align}
As the SDE and ODE approaches have the same marginal distributions, they are equivalent in probability when the score $\nabla_{\bsm(t)} \log \pi_{t}(\bsm(t))$ is known. However, the sample trajectories realized by the two formulations differ, and both empirical and theoretical results have provided evidence that the SDE based approach is more robust to score approximation error \cite{lu2022maximum, song2021maximum}.

Both the reverse SDE and probability flow ODE depend on $\nabla_{\bsm(t)} \log \pi_{t}(\bsm(t))$, the score (the gradient of the log-density) of the time-$t$ marginal distribution of $\bsm(t)$. These vector fields are \textit{a priori} unknown but can be learned from the provided samples using a parameterized model $s_{\bstheta}(\bsm(t), t): \mathbb{R}^D \times \mathbb{R}_+ \to \mathbb{R}^D$ (the \textit{score model}) with parameters $\bstheta \in \mathbb{R}^P$. We use $\bstheta^*$ to denote the optimized parameters of the model, which are obtained by minimizing the following objective over the set of samples $\mathcal{M}$ \cite{song2020score}:
\begin{equation}
L(\bstheta) = \frac{1}{2} \int_{T_{\mathrm{min}}}^T w(t) \sum_{i=1}^{N_s} \; \mathbb{E}_{\bsz} \; \left \lVert s_{\bstheta}(\bsm_i(t); \sigma(t)) + \frac{\bsz}{\sigma(t)} \right \rVert_2^2 \; dt.
\label{eq:denoising_score_matching}
\end{equation}
In the above objective, which is based on the denoising score matching objective originally introduced by Vincent \cite{vincent2011connection}, $w(t)$ is a specified weighting function, $\bsm_i(t) = \bsm_i + \sigma(t) \bsz $ is a sample from the noising distribution $\pi_{t \mid 0} (\bsm(t) \mid \bsm(0) = \bsm_i)$, $\mathbf{z}$ is white noise, and the lower integral limit $T_{\mathrm{min}} \geq 0$ is needed to ensure the objective is well-conditioned. 

\subsection{HIJACKING APPROACHES}
\label{subsec:hijacking}

The ascendance of diffusion models as a state-of-the-art generative modeling technique has spurred significant research on their use in the context of Bayesian inverse problems. In particular, there has been substantial interest in leveraging diffusion models trained on the prior distribution of a given inverse problem to sample from the posterior distribution, i.e., setting $\pi_0 = \pi_{\mathrm{pr}}$; see \cite{darassurvey} for a survey of these approaches. This general strategy has the advantage of enabling the same diffusion model to be used to sample from many posterior distributions corresponding to different choices of likelihood function without retraining the score model. In the remainder of this section, we discuss one algorithmic framework, which we refer to as the hijacking class of approaches,\footnote{These approaches are referred to as ``Explicit Approximations for Measurement Matching'' approaches in \cite{darassurvey}.} that follows this general strategy.

The main idea of the hijacking approaches is as follows: Given a diffusion model for the prior distribution, replace the $\nabla_{\bsm(t)} \log \pi_{t}(\bsm(t))$ term in \eqref{eq:reverse_SDE} with $\nabla_{\bsm(t)} \log \pi_{t \mid \bsy}(\bsm(t) \mid \bsy)$ to sample from the posterior, and expand the $ \pi_{t \mid \bsy}(\bsm(t) \mid \bsy)$ term using Bayes' Theorem. This yields the following reverse SDE:

\begin{align}
d\bsm(t) &= -2 \dot{\sigma}(t)\sigma(t) \nabla_{\bsm(t)} [\log \pi_{t}(\bsm(t))  \nonumber \\
&\; + \log \pi_{\bsy \mid t}(\bsy \mid \bsm(t))] \; dt +  \sqrt{2 \dot{\sigma}(t)\sigma(t)} \; d\bsw(t). \label{eq:reverse_SDE_hijacking}
\end{align}
In the above equation, an estimate of $\nabla_{\bsm(t)} \log \pi_{t}(\bsm(t))$ is given by the score model of the diffusion model pre-trained on the prior distribution. The $\nabla_{\bsm(t)} \log \pi_{\bsy \mid t}(\bsy \mid \bsm(t))$ term, which is known as the \textit{noisy likelihood function} \cite{zhang2024improving}, then ``hijacks'' the pre-trained diffusion process with information provided by the measurements. However, while the distribution of $\bsy$ given $\bsm(0)$ is given by the likelihood function, the noisy likelihood function is not known. In particular, computation of the noisy likelihood function
\begin{align}
\label{eq:y_given_mt}
&\pi_{\bsy \mid t}(\bsy \mid \bsm(t)) \nonumber \\
& \quad = \int \pi_{\mathrm{like}} (\bsy \mid \bsm(0)) \, \pi_{0 \mid t} (\bsm(0) \mid \bsm(t)) \; d\bsm(0),
\end{align}
is in general intractable.
In practice this issue is resolved by using simple approximations of the \textit{denoising distribution} $\pi_{0 \mid t} (\bsm(0) \mid \bsm(t))$ \cite{chung2023diffusion, boys2023tweedie, song2023pseudoinverseguided, yu_2023_freedom}. For example, in \cite{chung2023diffusion}, the denoising distribution is modeled as a Dirac delta with mass centered on $\mathbb{E}[\bsm(0) \mid \bsm(t) ]$, and Tweedie's formula \cite{efron2011tweedie} is employed to compute the expectation. In the approach of Boys et al \cite{boys2023tweedie}, a Gaussian approximation is employed, with the mean and covariance of the Gaussian computed using a generalized version of Tweedie's formula \cite{meng2021estimating}. However, this approach is only applicable when the likelihood function is linear-Gaussian, as otherwise the integral in \eqref{eq:y_given_mt} cannot be straightforwardly computed. 

The approaches in \cite{chung2023diffusion} and \cite{boys2023tweedie} discussed above use simple approximations of the denoising distribution to make the integral in \eqref{eq:y_given_mt} tractable and enable approximation of the noisy likelihood function. In particular, both approaches use unimodal approximations. While this may be a good approximation for many prior distributions of interest when $t \approx 0$, if the prior is multimodal then $\pi_{0 \mid t} (\bsm(0) \mid \bsm(t))$ will be multimodal as well for large enough $t$. There can therefore be significant errors in the approximation of the noisy likelihood function for large $t$, which in turn induces errors in the distribution of $\bsm(t)$ in \eqref{eq:reverse_SDE_hijacking}. These errors are difficult for subsequent hijacking iterations to correct, as discretizations of the hijacking reverse SDE use small step sizes $\Delta t > 0$ that ensure $\bsm(t - \Delta t)$ will be close to $\bsm(t)$. In practice, it has been observed \cite{zhang2024improving} that this can lead to poor performance on certain inverse problems (in particular, nonlinear inverse problems), with samples obtained that are consistent with the likelihood function but lie in low-density regions with respect to the prior. These issues have motivated the development of a class of methods that address this issue by decoupling $\bsm(t)$ and $\bsm(t - \Delta t)$, which we now introduce.

\section{PROPOSED APPROACH}
\label{sec:proposed}

In this section, we introduce the proposed BIPSDA framework for solving Bayesian inverse problems using diffusion models. Unlike the hijacking approaches discussed in the previous section, our framework uses a recently introduced technique called \textit{decoupled noise annealing} \cite{zhu2023denoising, zhang2024improving} to avoid approximation of the noisy likelihood function and prevent samples from getting stuck in low-density regions. In what follows, we first introduce the framework, which generalizes previously introduced decoupled noise annealing approaches \cite{zhang2024improving, zhu2023denoising} and can yield new algorithms through flexible combinations of the framework's design choices. After discussing the relationship between our framework and previously proposed decoupled noise annealing approaches, we give examples of concrete algorithms that can be realized within our framework. 

\subsection{THE FRAMEWORK}

As in other decoupled noise annealing approaches \cite{zhu2023denoising, zhang2024improving}, the proposed BIPSDA framework generates a sequence of iterates that are approximate samples from the noisy posterior distribution at a series of decreasing noise levels. 
Specifically, each iteration is comprised of two stages: a \textit{prediction} stage that generates a posterior sample that is consistent with the current iterate and the measurements, and a \textit{corruption} stage in which noise is added back to the predicted posterior sample. 

Concretely, given the current iterate $\bsm(t)$ and the measurements $\bsy$, the goal of the prediction stage is to obtain an approximate sample $\bsm(0)$ from the \textit{prediction distribution} $\pi_{0 \mid t, \bsy} (\bsm(0) \mid \bsm(t), \bsy)$. In the corruption stage, $\bsm(t - \Delta t)$ is obtained by adding noise back to $\bsm(0)$, i.e., by sampling from the Gaussian noising distribution $\pi_{t - \Delta t \mid 0}(\bsm(t - \Delta t) \mid \bsm(0))$. In \cite{zhang2024improving}, Zhang et al prove that if $\bsm(t)$ is a sample from $\pi_{t \mid \bsy}(\bsm(t) \mid \bsy)$, then this procedure yields a sample from $\pi_{t - \Delta t \mid \bsy}(\bsm(t - \Delta t) \mid \bsy)$. For sufficiently large $T$, once can assume that $\pi_{T \mid \bsy} (\bsm(T) \mid \bsy) \approx \mathcal{N}(\bsm(T); \mathbf{0}, \sigma^2(T) \mathbf{I})$. So starting from $\bsm(T)$, the prediction and corruption steps can be iteratively applied with the timestep annealed from $ T$ down to $0$ to yield a sample from the posterior distribution. 

The prediction stage of the algorithm outlined above requires sampling from the prediction distribution $\pi_{0 \mid t, \bsy} (\bsm(0) \mid \bsm(t), \bsy)$. Using Bayes' rule and the conditional independence of $\bsy$ from $\bsm(t)$ given $\bsm(0)$, it holds that \cite{zhang2024improving}: 
$$
\pi_{0 \mid t, \bsy} (\bsm(0) \mid \bsm(t), \bsy) \propto \pi_{\mathrm{like}}(\bsy \mid \bsm(0)) \; \pi_{0 \mid t}(\bsm(0) \mid \bsm(t)). 
$$
Unfortunately, while the likelihood function is known, the denoising distribution $\pi_{0 \mid t}(\bsm(0) \mid \bsm(t))$ is not and requires approximation. The prediction stage can therefore be broken down into two substages: approximation of the denoising distribution and sampling from the prediction distribution using the approximate denoising distribution. Here it is worth noting that the hijacking approaches discussed in the previous section also require approximation of the denoising distribution, and approximations used in the hijacking setting can also be used in the BIPSDA setting (e.g., the Tweedie formula based Gaussian approximation of Boys et al \cite{boys2023tweedie}). However, approximations in the hijacking setting must be chosen so that the integral in \eqref{eq:y_given_mt} and the subsequent score operation in \eqref{eq:reverse_SDE_hijacking} can be evaluated efficiently. In our setting, there are no such restrictions. In what follows, we use $\pi_{\mathrm{aprx}}(\bsm(0) \mid \bsm(t))$ to denote an approximation of the denoising distribution $\pi_{0 \mid t}(\bsm(0) \mid \bsm(t))$. 

In summary, the BIPSDA framework requires three steps at each iteration: approximation of the denoising distribution, sampling from the prediction distribution using the approximate denoising distribution, and corruption of the predicted sample using the Gaussian noising distribution. A full outline is provided in Algorithm \ref{alg:bipsda}. In what follows, we discuss previously proposed methods in the literature that fall within the BIPSDA framework, with particular attention paid to the techniques used for approximating the denoising distribution and for sampling from the approximate prediction distribution. We then give examples of novel algorithms that can be realized in our framework through different choices of the approximation distribution and sampling scheme.

\begin{algorithm}[t]
	\caption{ Bayesian Inverse Problem Solvers through Diffusion Annealing (BIPSDA)}
	\label{alg:bipsda}
	\begin{algorithmic}[1]
        \vspace{1mm}
		\STATE \textbf{Input:} Decreasing timesteps $[t_{N_A}, t_{N_A - 1}, \cdots, t_0]$, likelihood function $\pi_{\mathrm{like}}$, noise schedule $\sigma(t)$, trained score model $s_{\bstheta^*}(\bsm(t), t)$ \vspace{1mm}
		\STATE \textbf{Output:} Approximate sample from the posterior distribution $\pi_{\mathrm{post}}(\bsm \mid \bsy)$ \\ \vspace{1mm}
		\STATE Initialize $\bsm(t_{N_A}) \sim \mathcal{N}(\bsm(t_{N_A}); \mathbf{0}, \sigma^2(t_{N_A})\mathbf{I})$  \vspace{1mm}
		\FOR {$i = N_A, N_A - 1, \dots, 1$}\vspace{1mm}
		\STATE Compute $\pi_{\mathrm{aprx}}$ as approximation of $\pi_{0 \mid t}$ obtained using score model  \vspace{1mm}
		\STATE Sample $\bsm(0) \sim \pi_{\mathrm{like}}(\bsy \mid \bsm(0)) \mkern 1mu \pi_{\mathrm{aprx}}(\bsm(0) \mid \bsm(t_i))$ \vspace{-1mm}
        \STATE Sample $\bsm(t_{i-1}) \sim \mathcal{N}(\bsm(t_{i-1}); \bsm(0), \sigma^2(t_{i-1}) \mathbf{I})$ \vspace{1mm}
		\ENDFOR \vspace{1mm}
        \STATE Return $\bsm(t_0)$ \vspace{1mm}
	\end{algorithmic}
\end{algorithm}

\subsection{RELATIONSHIP TO PRIOR WORK}

Over the last two years, several diffusion-based Bayesian inverse problem solvers have been proposed that use decoupled noise annealing. Of these approaches, the work of Zhang et al \cite{zhang2024improving} is the most closely related to the proposed BIPSDA framework. In their approach, referred to as Decoupled Annealing Posterior Sampling (DAPS), the approximation of the denoising distribution  (line 5 in Algorithm \ref{alg:bipsda}) takes the form of a Gaussian distribution with mean $\bsm_{\mathrm{aprx}}$ and covariance $\mathbf{C}_{\mathrm{aprx}} = \beta(t)^2 \mathbf{I}$, where the noise level $\beta(t)$ is chosen using heuristics. 
The mean of $\bsm_{\mathrm{aprx}}$ is set as an estimate of $\mathbb{E}_{0 \mid t}[\bsm(0) \mid \bsm(t)]$, with the estimate obtained by solving the probability flow ODE in \eqref{eq:prob_flow_ODE} with initial value $\bsm(t)$. 
To sample from the corresponding approximate prediction distribution (line 6 in Algorithm \ref{alg:bipsda}), the DAPS approach uses MCMC algorithms such as the Euler–Maruyama method for discretizing Langevin dynamics (without Metropolis adjustment), or the Hamiltonian Monte Carlo algorithm \cite{neal2012mcmc}.

Another approach that is closely related to the present work is the DiffPIR algorithm of Zhu et al \cite{zhu2023denoising}. Like the DAPS algorithm, this approach models the denoising distribution as a Gaussian, with the mean set as an estimate of the conditional mean $\mathbb{E}_{0 \mid t}[\bsm(0) \mid \bsm(t)]$ and the covariance matrix chosen to be a scalar multiple of the identity.  However, unlike DAPS, DiffPIR estimates the conditional mean  using Tweedie's formula, which is exact up to error in the score model. 
Another key difference between DAPS and DiffPIR is that DiffPIR
 does not use Langevin dynamics to sample from the corresponding prediction distribution. Instead, inspired by plug-and-play algorithms for image restoration (see, e.g., \cite{chan2016plug}),
 it solves a
 maximum a posteriori (MAP) estimation problem corresponding to the prediction distribution.  This problem can be efficiently solved using fast numerical optimization methods or, in some cases, close-form proximal operators \cite{parikh2014proximal}. While lacking theoretical guarantees, it is empirically observed that DiffPIR can produce very diverse samples (see Figure 5 in \cite{zhu2023denoising}). 

Other diffusion-based Bayesian inverse problem solvers that use decoupled noise annealing include those proposed in \cite{alkhouri2024sitcom} and \cite{song2023solving}. In \cite{alkhouri2024sitcom}, a prediction-corruption decoupled noise annealing approach, dubbed SITCOM, was proposed. As in the proposed BIPSDA framework, the corruption stage of SITCOM is implemented by sampling from the Gaussian noising distribution. However, the prediction stage differs. In BIPSDA, the prediction stage aims to sample from the prediction distribution $\pi_{0 \mid t, \bsy} (\bsm(0) \mid \bsm(t), \bsy)$. In the SITCOM approach, the prediction stage requires first evaluating a proximal operator to obtain a point $\bsm(t)$ that is consistent with the measurement $\bsy$ and close to the current iterate $\bsm(t)$. Tweedie's formula is then applied to $\bsm(t)$ to obtain the prediction $\bsm(0)$. In \cite{song2023solving}, an approach similar to the DiffPIR algorithm was proposed. Unlike DiffPIR, however, which solves for the MAP point of the prediction distribution, \cite{song2023solving} proposes the use an iterative algorithm to maximize the likelihood function, with $\bsm_{\mathrm{aprx}}$ used to initialize the iterative algorithm.

\subsection{SPECIFIC VARIANTS OF THE FRAMEWORK}
\label{subsec:example_algs}

The proposed BIPSDA framework contains many novel algorithms for solving Bayesian inverse problems with diffusion models that can be realized through different choices of the approximation distribution and sampling scheme. In the following, we discuss these design choices and provide an outline of the algorithms that will be tested in the numerical studies. 
\subsubsection{Approximation of the denoising distribution} In this work, we consider three different approaches to build a Gaussian approximation to the denoising distribution (c.f. line 5 in Algorithm \ref{alg:bipsda}). That is, we consider approximations to $\pi_{0 \mid t}(\bsm(0) | \bsm(t))$ of the form
\begin{equation}
\label{eq:gaussian_denoising_approx_III}
\pi_{\mathrm{aprx}}(\bsm(0) \mid \bsm(t)) = \mathcal{N}(\bsm(0); \bsm_{\mathrm{aprx}}, \mathbf{C}_{\mathrm{aprx}}).
\end{equation}
The rationale for restricting our focus to a Gaussian approximation is that,  while in principle arbitrarily complex approximations of the denoising distribution can be realized by first sampling from the denoising distribution $\pi_{0 \mid t}(\bsm(0) \mid \bsm(t))$ using the reverse SDE and then fitting a model to the samples, these approximations may be computationally prohibitive in practice. The approaches considered are denoted as `ODE', Tweedie Uncorrelated (`TU'), and Tweedie Correlated (`TC').

In the `ODE' approach, which is based on the DAPS algorithm, the mean $\bsm_{\mathrm{aprx}}$ is obtained by solving the probability flow ODE in Eq. \eqref{eq:prob_flow_ODE}, and the covariance is chosen as $\mathbf{C}_{\mathrm{aprx}} = \beta(t)^2 \mathbf{I}$, with $\beta(t) = \mathcal{O}(\sigma(t))$. User-provided hyperparameters of this approach are the marginal variance of the denoising distribution $\beta(t)^2$ and the schedule of time steps in the discretization of the probability flow ODE.

In the `TU' approach, which is based on the  DiffPIR algorithm, the mean is obtained using Tweedie's formula, i.e., 
$$
\bsm_{\mathrm{aprx}} = \bsm(t) + \sigma^2(t)  s_{\bstheta^*}(\bsm(t), t),
$$
and the covariance is set as $\mathbf{C}_{\mathrm{aprx}} = \beta(t)^2 \mathbf{I}$, with $\beta(t) = \mathcal{O}(\sigma(t))$ the only user-provided hyperparameter. We note that, as shown in Appendix \ref{sec:appendix}, under a particular choice of parametrization and discretization, the probability flow ODE estimate used in the `ODE' approach can be made to coincide with Tweedie formula.

Finally, in the `TC' approach, which borrows some ideas from the hijacking approach of Boys et al \cite{boys2023tweedie}, the mean $\bsm_{\mathrm{aprx}}$ is set as in the `TU' approach, while the covariance is set using the generalized version of Tweedie's formula \cite{meng2021estimating}, i.e., 
$$
\mathbf{C}_{\mathrm{aprx}} = \sigma^2(t) \left[ \mathbf{I} + \sigma^2(t) \nabla_{\bsm(t)} s_{\bstheta^*}(\bsm(t), t) \right]. 
$$
This variant is free from user-provided hyperparameters, but requires access to accurate estimates of the Jacobian of the noisy prior score.

\subsubsection{Sampling from the prediction distribution}
The BIPSDA framework also provides considerable flexibility regarding the choice of the sampling scheme. Three variants, denoted as Langevin (`Lang'), maximum a posteriori estimation (`MAP'), and randomize-then-optimize (`RTO'), are considered here.

The `Lang' variant is based on  DAPS \cite{zhang2024improving} and employs Markov chain Monte Carlo (MCMC) algorithms, such as Langevin dynamics or Hamiltonian Monte Carlo \cite{neal2012mcmc}, to sample the prediction distribution. When using uncorrected Langevin dynamics the step-size and the number of time steps are key user-provided hyperparameters that control the trade-off between accuracy of the samples and computational costs.

The `MAP' variant is based on DiffPIR \cite{zhu2023denoising} and simply computes the maximum a posteriori (MAP) point of the prediction distribution using a deterministic optimization algorithm, rather than drawing samples. While lacking solid theoretical foundation, this approach is computationally efficient. User-provided hyperparameters depends on the specific choice of the optimizer used, and, as long as the optimization problem is solved with sufficient accuracy, have limited effect on the quality of the solution.

Finally, the `RTO' variant is a novel contribution of our work. It is inspired by the randomize-then-optimize (RTO) \cite{bardsley2014randomize,ba2022randomized,Kainan2018randomized} algorithm for approximate sampling from the prediction distribution. Here the key idea is to obtain an approximate sample from the prediction distribution by first adding noise to both the measurement and the denoising distribution mean, then solving for the MAP point of the corresponding noise-perturbed prediction distribution. Concretely, assuming the Gaussian approximation of the denoising distribution in Eq. \eqref{eq:gaussian_denoising_approx_III} and the Gaussian additive noise model in Eq. \eqref{eq:additive-gaussian}, 
`RTO' generates an approximate sample $\bsm(0)$ from the prediction distribution by solving
\begin{equation}
\bsm(0) = \operatornamewithlimits{argmax}_{\bsm} \pi_{\rm like}( \bsy' | \bsm ) \; \mathcal{N}(\bsm; \bsm_{\mathrm{aprx}}', \mathbf{C}_{\mathrm{aprx}}),
\end{equation}
where 
\begin{align*}
\bsm_{\mathrm{aprx}}' &\sim \mathcal{N}(\bsm_{\mathrm{aprx}}'; \bsm_{\mathrm{aprx}}, \mathbf{C}_{\mathrm{aprx}}), \text{ and}\\
\bsy' & \sim \mathcal{N}(\bsy'; \bsy, \boldsymbol{\Sigma}_{\bsz})
\end{align*}
denote the noise-perturbed mean and measurement, respectively.
Note that this approach can be viewed as an application of Algorithm 1 in \cite{Kainan2018randomized} to the subproblem of sampling from the prediction distribution (line 6 in Algorithm \ref{alg:bipsda}). Like the `MAP' variant used in DiffPIR, it has the advantage of enabling efficient MAP solvers to be employed for sampling from prediction distribution. Further, in the linear-Gaussian likelihood setting, where the corresponding prediction distribution is also Gaussian, it can be proven that this approach corresponds to exact sampling from the prediction distribution \cite{bardsley2014randomize,bardsley2012mcmc}.

\subsubsection{Analyzed variants}

Mixing and matching different variants to approximate the denoising distribution with those to sample from the prediction distribution, we obtain the 9 variants summarized in Table \ref{tab:example_algs}. The DAPS and DiffPIR algorithms correspond to the `Lang-ODE' and `MAP-TU' variants of our framework, respectively.

\begin{table}
\centering
    \begin{tabular}{c c|ccc}
    & &  \multicolumn{3}{c}{Denoising Dist. Approximation} \\ 
 & & \hspace{.4cm} ODE \hspace{.4cm} & \hspace{.4cm} TU \hspace{.4cm} & \hspace{.4cm} TC \hspace{.4cm} \\ \hline
 \parbox[t]{2mm}{\multirow{3}{*}{\rotatebox[origin=c]{90}{Sampling Alg.}}} \hspace{.15cm} & Lang \hspace{.15cm} & Lang-ODE & Lang-TU & Lang-TC \\
   & MAP \hspace{.15cm}& MAP-ODE & MAP-TU & MAP-TC\\
   & RTO \hspace{.15cm} & RTO-ODE & RTO-TU & RTO-TC \\
    \end{tabular}
    \vspace{0.15in}
    \caption{Visual representation of the nine example algorithms tested, each of which corresponds to different choices of the denoising distribution approximation (line 5 in \ref{alg:bipsda}) and sampling scheme (line 6 in \ref{alg:bipsda}). Note that `Lang-ODE' is the DAPS algorithm, while `MAP-TU' is the DiffPIR algorithm; the other seven variants are novel.}
    \label{tab:example_algs}
\end{table}

\section{NUMERICAL STUDIES}
\label{sec:numerical_staples}

The performance of diffusion-model-based posterior sampling algorithms is often evaluated on high-dimensional problems (e.g., imaging problems), where the prior distribution is characterized implicitly through a dataset of examples (samples). In this setting, the posterior density is analytically unknown, making it difficult to evaluate the ability of the algorithm under consideration to capture the global structure of the posterior distribution and facilitate uncertainty quantification. To address this issue, we created a suite of benchmark problems for which the prior distribution, as well as the score of the denoising distribution, is analytically known. Using these benchmark problems, which use a Gaussian mixture as the prior distribution, we can rigorously assess the performance of the algorithms in our framework, which include the popular DAPS \cite{zhang2024improving} and DiffPIR \cite{zhu2023denoising} algorithms. 

The use of a Gaussian mixture prior to test the performance of diffusion-model-based posterior samplers has a number of advantages. First, in this setting the posterior density function is known, and in certain cases (e.g., linear-Gaussian likelihood function, where the posterior is also a Gaussian mixture) can be straightforwardly sampled from, enabling us to investigate errors in the global structure of the samples from these algorithms. Second, under this choice $\pi_{t}(\bsm(t))$ 
 is also a Gaussian mixture and the noisy prior score $\nabla_{\bsm(t)} \log \pi_{t}(\bsm(t))$ is known analytically. This enables the decoupling of error in the score modeling error from error inherent to the sampling algorithm under consideration in a principled manner. Third, Gaussian mixture distributions can be made arbitrarily complex through the choice of the number of modes and the covariances of each component Gaussian.

Each of the studies conducted corresponds to a different choice of likelihood function for the Bayesian inverse problem. In the first and second studies, the likelihood function was inspired by image inpainting problems. We assumed a (linear) binary sampling mask as the forward model and additive Gaussian measurement noise. In the first study a low noise regime was considered, with the noise level set so that the posterior was unimodal up to numerical precision. In the second study, a high noise regime was considered where the posterior contained distinct modes. The third study considered a nonlinear likelihood function with Poisson measurement noise inspired by problems in X-ray computed tomography (CT) \cite{tau2020proximal}. Finally, the fourth study considered a problem with a nonlinear Gaussian phase retrieval forward model and Gaussian measurement noise. Phase retrieval is a difficult inverse problem often arising in optical imaging applications \cite{shechtman2015phase} and exhibits a complicated multi-modal posterior distribution.

In each study, we tested the performance of the nine algorithms within the BIPSDA framework summarized in Table \ref{alg:bipsda}.
Additional details regarding the implementation of these algorithms, as well as the problem setting and metrics used to evaluate the performance of the algorithms, are provided in the following subsections. 

\subsection{PROBLEM SETTINGS}

\begin{figure}
    \centering
    \includegraphics[trim = {0, 0, 0, 2.65cm}, clip, width=0.8\linewidth]{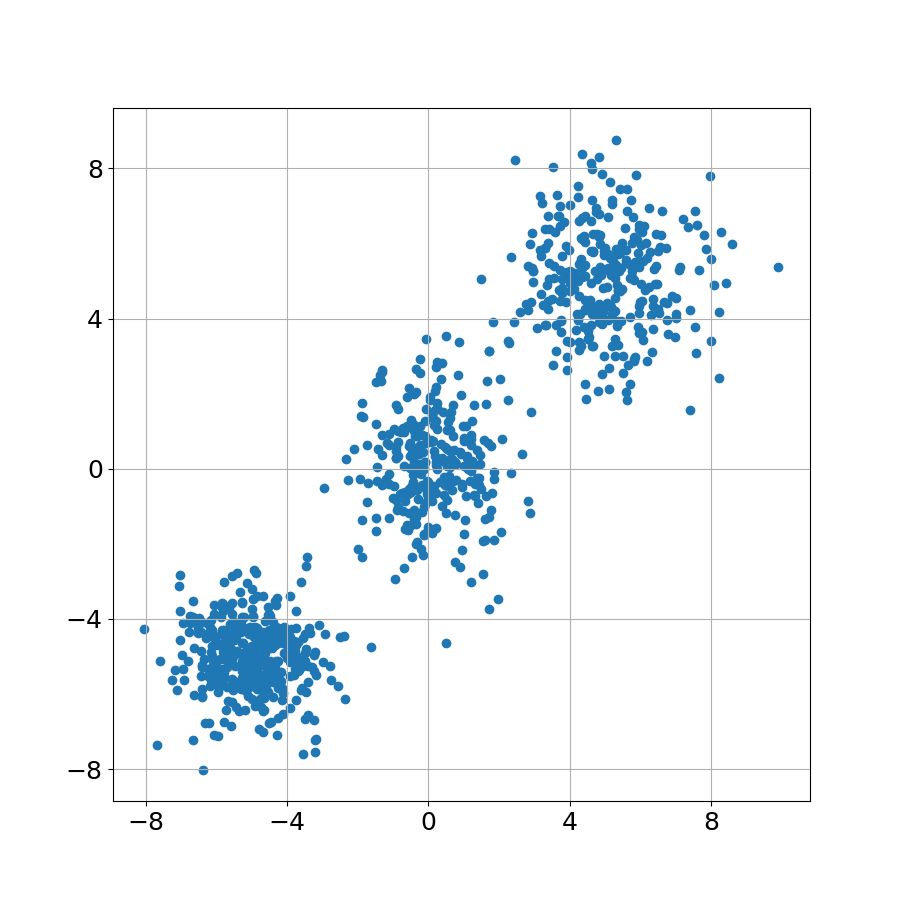}
    \caption{Two components of samples from the ten-dimensional Gaussian mixture prior distribution used in all of the numerical studies. The distribution has three distinct modes and each mode has a different covariance matrix. }
    \label{fig:prior-plot}
\end{figure}

As previously mentioned, all of the numerical studies described in this work consider inverse problems with inversion parameter $\bsm \in \mathbb{R}^{10}$ and a Gaussian mixture prior, i.e., 
$$
\pi_{\mathrm{pr}}(\bsm) = \sum_{i=1}^{N_m} w_i \, \pi_{\mathrm{pr}, i}(\bsm),
$$
where each $\pi_{\mathrm{pr}, i}$ is a Gaussian, $w_i$ is the weight of the $i$th component, and we set $N_m = 3$. 
Here $\pi_{\mathrm{pr}, 1}$ has mean $[-5, \cdots, -5]^T$ and identity covariance matrix, $\pi_{\mathrm{pr}, 2}$ has mean $[0, \cdots, 0]^T$ and a diagonal covariance matrix with entries linearly spaced between $1$ and $2$, and $\pi_{\mathrm{pr}, 3}$ has mean $[5, \cdots, 5]^T$ and a covariance matrix with the same eigenvalues as $\pi_{\mathrm{pr}, 2}$ but randomly chosen eigenvectors. The corresponding weights were set as $w_1=0.4$, $w_2 = 0.3$, and $w_3 = 0.3$. A plot of two components of samples from this distribution is shown in Figure \ref{fig:prior-plot}. 

Both the stylized inpainting and phase retrieval studies consider likelihood functions with additive Gaussian noise, as in \eqref{eq:additive-gaussian}, with white noise covariance $\boldsymbol{\Sigma}_{\mathbf{z}} = \tau^2 \mathbf{I}$. In the stylized inpainting studies, the forward model was given as $f(\bsm) = \mathbf{A} \bsm$, where $\mathbf{A} \in \mathbb{R}^{8 \times 10}$ is a binary subsampling operator. The noise standard deviation was set as $\tau = 0.1$ (low noise regime, $\mathrm{SNR} \approx 30.4687 \; \mathrm{dB}$) or $\tau = 5$ (high noise regime,  $\mathrm{SNR} \approx -3.5107 \; \mathrm{dB}$). For the stylized phase retrieval problem, the forward model was given as $f(\bsm) = (\mathbf{B} \bsm)^2$, where here the square operation is applied elementwise and $\mathbf{B} \in \mathbb{R}^{5 \times 10}$ is a matrix with i.i.d. standard Gaussian entries. The noise level for this problem was set as $\tau = 25$ ($\mathrm{SNR} \approx -0.8243 \; \mathrm{dB}$). Finally, the likelihood function of the stylized x-ray tomography has the form
$$
\pi_{\mathrm{like}} (\bsy \mid \bsm) = \mathrm{Poi}(f(\bsm)), \quad f(\bsm) = I_0 \; \text{with } \mathrm{exp} \left ( - \mathbf{C} \mathbf{m} \right),
$$
where $\mathrm{Poi}(\cdot)$ is the Poisson distribution, the exponential operator is applied pointwise, and $\mathbf{C} \in \mathbb{R}^{15 \times 10}$. We set $I_0 = 1000$ and populated $\mathbf{C}$ with i.i.d. random entries from the uniform distribution over the interval $[.01, .05]$ ($\mathrm{SNR} \approx 23.1709 \; \mathrm{dB}$).

\subsection{ALGORITHM IMPLEMENTATIONS}
\label{subsec_implementation}

Each of the BIPSDA algorithms variants summarized in Table \ref{tab:example_algs} was implemented using $N_{A} = 200$ timepoints in the noise annealing, with the $t_i$ set by polynomial interpolation between $T$ and $0$. The `ODE' and `TU' variants for approximating the denoising distribution were implemented with $\beta(t) = \sigma(t)$ in the covariance approximation $\mathbf{C}_{\mathrm{aprx}} = \beta(t)^2 \mathbf{I}$; the `ODE' method for approximating the denoising distribution was implemented by solving the probability flow ODE using the Euler discretization with five discretization time-steps. The above settings are the same as used by Zhang et al; see \cite{zhang2024improving} for details. 

The implementation of the Langevin dynamics based sampler and the MAP solver used in the `MAP` and `RTO` variants was problem specific. The sampler for the stylized inpainting problem was implemented with step size $= 5 \times 10^{-5}$, $100$ subiterations per algorithm iteration, and no Metropolis correction, as in \cite{zhang2024improving}. The implementation was the same for the stylized x-ray tomography based problem, but with $1000$ subiterations instead of $100$. For the phase retrieval problem, we observed severe stability issues with unadjusted Langevin dynamics, and instead implemented Langevin dynamics with Metropolis adjustment and preconditioning \cite{girolami2011riemann}. Here $1000$ subiterations were used with step size $=.2$, and at each iteration we fixed the preconditioning matrix as the average of a set of matrices, where each matrix in the set is the Gauss-Newton Hessian approximation \cite{nocedal1999numerical} evaluated at different representative sample points. Regarding the MAP solver used in the `MAP` and `RTO` variants, for the stylized inpainting problem we exploited the fact that the MAP point has a closed-form expression. For the stylized x-ray tomography and phase retrieval problems, we estimated the MAP point using the PyTorch implementation of the limited-memory Broyden–Fletcher–Goldfarb–Shanno (L-BFGS) algorithm \cite{nocedal1999numerical} (40 subiterations per algorithm iteration, strong Wolfe line search).

Each of the algorithms described above was implemented using the same pretrained score model corresponding to a diffusion model with $\sigma(t) = t$ and $T = 10$. The score model was given as $s_{\bstheta}(\bsm(t), t) = \nabla_{\bsm(t)} g_{\bstheta}(\bsm(t), t)$ where $g_{\bstheta}(\bsm(t), t): \mathbb{R}^{10} \times \mathbb{R}_{+} \to \mathbb{R}$ was a deep neural network. In particular, the network architecture had six total hidden layers (width $= 512$), with the noise level $\sigma(t)$ appended to the input $\bsm(t)$ and to the network state after the 4th layer, width $= 512$. The network weights were trained to minimize the loss function in \eqref{eq:denoising_score_matching}, with $N = 80,000$ total samples in the training set, $w(t) = \sigma^2(t)$, and $T_{\mathrm{min}} = .01$. The Adam optimizer \cite{kingma2014adam} was used with learning rate $10^{-5}$ and a batch size of $8000$ $\{ \bsm_i, \bsz_i, t_i \}$ pairs (with the $\bsz_i$ and $t_i$ picked i.i.d. at each iteration). The network was trained for $50,000$ optimization iterations. Additionally, each algorithm was also tested with the exact ground truth score $\nabla_{\bsm(t)} \log \pi_{t} \bsm(t)$ used instead of the score model to isolate the effect of the score modeling error from other sources of error in the algorithms' performance. 

\subsection{EVALUATION METHODS}

We evaluated the performance of the BIPSDA algorithms by comparing the samples obtained by the algorithms with those produced by a reference method (exact sampling in the linear inverse problem case, MCMC for the non-linear inverse problems). In particular, in the stylized inpainting studies, we leveraged the fact that the ground truth posterior is a Gaussian mixture to obtain exact samples from the posterior as a reference. In the stylized x-ray tomography and phase retrieval studies, we used the PyDream implementation \cite{shockley2018pydream} of the MT-DREAM (ZS) algorithm \cite{laloy2012high}, a state-of-the-art MCMC posterior sampling algorithm for Bayesian inverse problems where the prior density function is known analytically, to obtain approximate ground truth samples from the posterior. In the stylized x-ray tomography study, we implemented PyDream with $10$ chains and $200,000$ iterations per chain, and obtained the final posterior samples by discarding the first half of each chain as burn-in and randomly sub-sampling the remaining samples. A similar procedure was followed in the stylized phase retrieval study, but with additional features incorporated to address the challenging nature of the problem. In particular, we took advantage of the fact that the posterior is a mixture distribution, i.e., 
\begin{align}
    \label{eq:post_decomp_pydream}
    \pi_{\mathrm{post}} (\bsm \mid \bsy) &\propto \pi_{\mathrm{like}} (\bsy \mid \bsm) \, \pi_{\mathrm{pr}}(\bsm) \nonumber \\
    &= \sum_{i=1}^{N_m} w_i \, \pi_{\mathrm{like}} (\bsy \mid \bsm) \, \pi_{\mathrm{pr}, i}(\bsm),
\end{align}
and sampled each of the $N_m$ components separately. For each component, we implemented PyDream with $40$ chains, $200,000$ iterations run per chain, and Latin hypercube sampling \cite{loh1996latin} to seed the sampling history. After discarding the first half of each chain as burn-in, the remaining samples were sub-sampled according to the weights of each of the $N_m$ components.

We assessed the convergence of the MCMC samples to the target distribution using  well-established MCMC diagnostics, including the potential scale reduction factor (PSRF) \cite{vats2021revisiting} and effective sample size (ESS) \cite{vehtari2021rank} metrics.  These metrics were computed for a randomly selected single posterior sampling trial. For the x-ray tomography problem, the maximum PSRF value over all parameter dimensions was 1.0006, with PSRF values less than 1.01 generally indicating that the parallel chains have stabilized, and samples are likely to have reached the target distribution \cite{vats2021revisiting}. The minimum ESS computed for each parameter independently was $47,077.65$, which is much larger than the number of posterior samples used in our evaluation metrics ($10,000$). For the phase retrieval study,  PSRF and ESS were computed separately for each PyDREAM run corresponding to each component of the posterior mixture in Eq. \eqref{eq:post_decomp_pydream}. The maximum PSRF value was $1.0004$ and the minimum ESS was $49,862.49$, which again indicate that a sufficient number of high-quality samples from the target distribution were generated.

The reference samples were compared to the samples produced by the BIPSDA algorithms both qualitatively and quantitatively. In particular, in each study we compared the performance of the BIPSDA algorithms to the reference over $100$ total trials, with each trial corresponding to sampling from the posterior distribution $\pi_{\mathrm{post}}(\bsm \mid \bsy_i)$ and $\bsy_i$ ($i=1, \ldots, 100$) chosen i.i.d. from the measurement distribution. In each trial, $10,000$ posterior samples were obtained from both the reference and each of the proposed algorithms. 

For the quantitative tests, four different error metrics were used: the central moment discrepancy (CMD) metric \cite{zellinger2017central}, the maximum mean discrepancy (MMD) metric \cite{gretton2012kernel}, and the two-norm error in both the predicted posterior mean and the predicted posterior pointwise variance. The CMD metric is the weighted sum of discrepancies between the central moments of the two distributions. In particular, for two distributions $\pi_1$ and $\pi_2$, the CMD metric can be written as follows: 
\begin{align}
\mathrm{CMD}(\pi_1, \pi_2) &= \frac{1}{\alpha} \left \lVert \mathbb{E}_{\pi_1}[\bsm] - \mathbb{E}_{\pi_2}[\bsm] \right \rVert_2  \\
& \quad + \sum_{k=2}^\infty \frac{1}{\alpha^k} \left \lVert c_k(\pi_1) - c_k(\pi_2) \right \rVert_2,\nonumber
\end{align}
where $c_k(\cdot)$ is the $k$th central moment of the given distribution and $\alpha > 0$ is a decay rate parameter. In practice, the infinite sum is truncated at some finite index $K$ (we use $K = 5$ in this work) for computability, which is theoretically justified by the fact that the terms in the sum can be shown to converge to zero as $k \to \infty$ for large enough $\alpha$. The expectations and central moments are also replaced by empirical approximations. Finally, in this work we set the decay rate as $\alpha = 4 \hat{\eta}_{\mathrm{max}}$, where $\hat{\eta}_{\mathrm{max}}$ is an empirical estimate of 
$$
\eta_{\mathrm{max}} = \mathbb{E}_{\bsy} \left [ 
|| \boldsymbol{\eta}(\bsy) ||_{\infty} \right]
$$
and $\boldsymbol{\eta}(\bsy) \in \mathbb{R}^D$ is the componentwise standard deviation of the posterior distribution $\pi_{\mathrm{post}} (\bsm \mid \bsy)$. 

The MMD metric is based on the difference between the two probability distributions in a reproducing kernel Hilbert space (RKHS) \cite{gretton2012kernel}. In particular, the MMD metric is an integral probability metric in the RKHS space and can be written as 
\begin{equation}
\mathrm{MMD}(\pi_1, \pi_2) = \sup_{f \in \mathcal{F}} \left(\mathbb{E}_{\pi_1}[f(\bsm)] - \mathbb{E}_{\pi_2}[f(\bsm)] \right),
\end{equation}
where $\pi_1$ and $\pi_2$ are probability distributions and $\mathcal{F}$ is the set of functions lying on the unit ball in the RKHS. An RKHS space is characterized by the choice of kernel function $k(\cdot, \cdot)$, with $ \{ k(\cdot, \bsm) \mid \bsm \in \mathbb{R}^D \}$ forming a set of basis functions for $\mathcal{F}$. In this work, we set the kernel as the sum of Gaussian radial basis functions with different bandwidths, i.e., 
$$
k(\bsm_1, \bsm_2) = \sum_{i=1}^{N_{b}} \mathrm{exp}\left( - \frac{|| \bsm_1 - \bsm_2 ||_2^2}{\epsilon_i} \right), 
$$
with $N_b = 5$, $\epsilon_i = \bar{\epsilon} \; 2^{i - \lceil N_b /2 \rceil } $, and $\bar{\epsilon}$ set as the average squared two-norm distance between samples from the reference distribution (as in \cite{gretton2012kernel}). 

\begin{table*}
\centering
    \begin{tabular}{c|cccc| cccc}
    \rowcolor{my_gray} &  \multicolumn{4}{c |}{Analytic Score} &  \multicolumn{4}{c}{Learned Score}  \\ 
    \rowcolor{my_gray} & \colorbox{my_gray}{\makecell{Mean \\ Error}}  & \colorbox{my_gray}{\makecell{Variance \\ Error}}  &  CMD   & MMD &  \colorbox{my_gray}{\makecell{Mean \\ Error}}  & \colorbox{my_gray}{\makecell{Variance \\ Error}}  &  CMD   & MMD  \\ \hline

Lang-TU &\makecell{$0.029$ \\ $(0.011, 0.055)$} & \makecell{$0.64$ \\ $(0.44, 0.86)$} & \makecell{$0.04$ \\ $(0.03, 0.06)$} & \makecell{$0.043$ \\ $(0.033, 0.054)$}&\makecell{$0.046$ \\ $(0.013, 0.092)$} & \makecell{$0.64$ \\ $(0.43, 0.87)$} & \makecell{$0.05$ \\ $(0.03, 0.07)$} & \makecell{$0.043$ \\ $(0.032, 0.055)$} \\ 
Lang-TC &\makecell{$0.024$ \\ $(0.009, 0.048)$} & \makecell{$0.64$ \\ $(0.43, 0.87)$} & \makecell{$0.04$ \\ $(0.03, 0.06)$} & \makecell{$0.043$ \\ $(0.033, 0.052)$}& -- -- & -- -- & -- -- & -- -- \\ 
Lang-ODE &\makecell{$0.032$ \\ $(0.016, 0.060)$} & \makecell{$0.10$ \\ $(0.02, 0.21)$} & \makecell{$\mathbf{0.01}$ \\ $(\mathbf{0.01}, \mathbf{0.02})$} & \makecell{$0.002$ \\ $(0.001, 0.003)$}&\makecell{$0.047$ \\ $(0.016, 0.090)$} & \makecell{$0.11$ \\ $(0.02, 0.25)$} & \makecell{$0.02$ \\ $(0.01, 0.03)$} & \makecell{$0.002$ \\ $(0.001, 0.005)$} \\ 
MAP-TU &\makecell{$\mathbf{0.018}$ \\ $(\mathbf{0.007}, \mathbf{0.031})$} & \makecell{$0.95$ \\ $(0.71, 1.19)$} & \makecell{$0.06$ \\ $(0.04, 0.08)$} & \makecell{$0.142$ \\ $(0.135, 0.149)$}&\makecell{$0.046$ \\ $(0.013, 0.093)$} & \makecell{$0.95$ \\ $(0.70, 1.21)$} & \makecell{$0.07$ \\ $(0.04, 0.09)$} & \makecell{$0.143$ \\ $(0.133, 0.155)$} \\ 
MAP-TC &\makecell{$\mathbf{0.018}$ \\ $(\mathbf{0.007}, \mathbf{0.030})$} & \makecell{$0.96$ \\ $(0.71, 1.20)$} & \makecell{$0.06$ \\ $(0.04, 0.08)$} & \makecell{$0.143$ \\ $(0.135, 0.151)$}& -- -- & -- -- & -- -- & -- -- \\ 
MAP-ODE &\makecell{$0.023$ \\ $(0.013, 0.037)$} & \makecell{$0.47$ \\ $(0.34, 0.61)$} & \makecell{$0.03$ \\ $(0.02, 0.04)$} & \makecell{$0.028$ \\ $(0.025, 0.032)$}&\makecell{$0.047$ \\ $(0.016, 0.090)$} & \makecell{$0.47$ \\ $(0.33, 0.64)$} & \makecell{$0.04$ \\ $(0.02, 0.06)$} & \makecell{$0.029$ \\ $(0.025, 0.035)$} \\ 
RTO-TU &\makecell{$0.022$ \\ $(0.008, 0.037)$} & \makecell{$0.05$ \\ $(0.02, 0.07)$} & \makecell{$\mathbf{0.01}$ \\ $(\mathbf{0.01}, \mathbf{0.01})$} & \makecell{$\mathbf{0.001}$ \\ $(\mathbf{0.000}, \mathbf{0.001})$}&\makecell{$\mathbf{0.044}$ \\ $(\mathbf{0.011}, \mathbf{0.089})$} & \makecell{$\mathbf{0.06}$ \\ $(\mathbf{0.02}, \mathbf{0.09})$} & \makecell{$\mathbf{0.01}$ \\ $(\mathbf{0.01}, \mathbf{0.02})$} & \makecell{$\mathbf{0.001}$ \\ $(\mathbf{0.001}, \mathbf{0.003})$} \\ 
RTO-TC &\makecell{$0.020$ \\ $(0.008, 0.031)$} & \makecell{$\mathbf{0.03}$ \\ $(\mathbf{0.01}, \mathbf{0.06})$} & \makecell{$\mathbf{0.01}$ \\ $(\mathbf{0.00}, \mathbf{0.01})$} & \makecell{$\mathbf{0.001}$ \\ $(\mathbf{0.000}, \mathbf{0.001})$}& -- -- & -- -- & -- -- & -- -- \\ 
RTO-ODE &\makecell{$0.027$ \\ $(0.013, 0.044)$} & \makecell{$1.05$ \\ $(0.79, 1.28)$} & \makecell{$0.08$ \\ $(0.05, 0.10)$} & \makecell{$0.051$ \\ $(0.047, 0.056)$}&\makecell{$0.046$ \\ $(0.016, 0.090)$} & \makecell{$1.01$ \\ $(0.80, 1.22)$} & \makecell{$0.08$ \\ $(0.06, 0.10)$} & \makecell{$0.050$ \\ $(0.043, 0.057)$} \\ 
Reference &\makecell{$0.020$ \\ $(0.009, 0.034)$} & \makecell{$0.03$ \\ $(0.01, 0.07)$} & \makecell{$0.01$ \\ $(0.00, 0.01)$} & \makecell{$0.001$ \\ $(0.000, 0.001)$}&\makecell{$0.020$ \\ $(0.009, 0.034)$} & \makecell{$0.03$ \\ $(0.01, 0.07)$} & \makecell{$0.01$ \\ $(0.00, 0.01)$} & \makecell{$0.001$ \\ $(0.000, 0.001)$} \\ 
    \end{tabular}
    \vspace{0.15in}
    \caption{Stylized inpainting study, low noise regime: Average errors in estimation of the mean and pointwise variance of the posterior distribution, and average central moment discrepancy (CMD) and maximum mean discrepancy (MMD) errors, for different BIPSDA methods (with both learned and analytic score). Interdecile ranges are also displayed in parentheses, and a reference that shows the average errors between two sets of i.i.d. samples from the ground truth posterior is included. The best-performing method with respect to each error measure is bolded. In the learned score case, with the exception of the  Tweedie-Correlated (`TC') approaches that were not implemented, all approaches performed well, with the `MAP' and `RTO' variants providing particularly strong performance.}
    \label{tab:inpaintinglown}
\end{table*}

\begin{table*}
\centering
    \begin{tabular}{c|cccc| cccc}
    \rowcolor{my_gray} &  \multicolumn{4}{c |}{Analytic Score} &  \multicolumn{4}{c}{Learned Score}  \\ 
    \rowcolor{my_gray} & \colorbox{my_gray}{\makecell{Mean \\ Error}}  & \colorbox{my_gray}{\makecell{ Variance \\ Error}}  &  CMD   & MMD & \colorbox{my_gray}{\makecell{Mean \\ Error}}  & \colorbox{my_gray}{\makecell{ Variance \\ Error}}  &  CMD   & MMD  \\ \hline

Lang-TU &\makecell{$10.594$ \\ $(1.918, 17.898)$} & \makecell{$44.17$ \\ $(33.88, 50.87)$} & \makecell{$3.14$ \\ $(1.93, 4.27)$} & \makecell{$0.389$ \\ $(0.188, 0.589)$}&\makecell{$10.076$ \\ $(1.613, 17.952)$} & \makecell{$43.99$ \\ $(33.20, 51.19)$} & \makecell{$3.05$ \\ $(1.93, 4.26)$} & \makecell{$0.380$ \\ $(0.145, 0.585)$} \\ 
Lang-TC &\makecell{$10.597$ \\ $(2.358, 18.025)$} & \makecell{$44.84$ \\ $(35.84, 51.02)$} & \makecell{$3.13$ \\ $(2.06, 4.28)$} & \makecell{$0.399$ \\ $(0.241, 0.593)$}& -- -- & -- -- & -- -- & -- -- \\ 
Lang-ODE &\makecell{$10.622$ \\ $(2.201, 17.217)$} & \makecell{$46.76$ \\ $(36.76, 53.64)$} & \makecell{$3.19$ \\ $(1.98, 4.24)$} & \makecell{$0.371$ \\ $(0.188, 0.506)$}&\makecell{$10.071$ \\ $(1.465, 17.538)$} & \makecell{$45.76$ \\ $(34.98, 53.04)$} & \makecell{$3.09$ \\ $(1.98, 4.26)$} & \makecell{$0.362$ \\ $(0.157, 0.514)$} \\ 
MAP-TU &\makecell{$\mathbf{0.246}$ \\ $(\mathbf{0.059}, \mathbf{0.492})$} & \makecell{$1.29$ \\ $(0.55, 2.05)$} & \makecell{$\mathbf{0.08}$ \\ $(\mathbf{0.04}, \mathbf{0.13})$} & \makecell{$0.130$ \\ $(0.088, 0.168)$}&\makecell{$\mathbf{0.262}$ \\ $(\mathbf{0.073}, \mathbf{0.569})$} & \makecell{$1.41$ \\ $(0.61, 2.21)$} & \makecell{$\mathbf{0.08}$ \\ $(\mathbf{0.04}, \mathbf{0.13})$} & \makecell{$0.128$ \\ $(0.085, 0.166)$} \\ 
MAP-TC &\makecell{$0.261$ \\ $(0.054, 0.664)$} & \makecell{$1.51$ \\ $(0.52, 2.61)$} & \makecell{$0.09$ \\ $(0.04, 0.17)$} & \makecell{$0.133$ \\ $(0.091, 0.164)$}& -- -- & -- -- & -- -- & -- -- \\ 
MAP-ODE &\makecell{$0.502$ \\ $(0.307, 0.719)$} & \makecell{$\mathbf{0.80}$ \\ $(\mathbf{0.34}, \mathbf{1.38})$} & \makecell{$0.10$ \\ $(0.06, 0.15)$} & \makecell{$0.035$ \\ $(0.022, 0.047)$}&\makecell{$0.479$ \\ $(0.324, 0.651)$} & \makecell{$\mathbf{0.99}$ \\ $(\mathbf{0.48}, \mathbf{1.93})$} & \makecell{$0.10$ \\ $(0.06, 0.16)$} & \makecell{$0.035$ \\ $(0.023, 0.047)$} \\ 
RTO-TU &\makecell{$1.156$ \\ $(0.270, 2.140)$} & \makecell{$5.04$ \\ $(1.20, 8.42)$} & \makecell{$0.34$ \\ $(0.12, 0.52)$} & \makecell{$0.016$ \\ $(0.002, 0.032)$}&\makecell{$1.017$ \\ $(0.198, 2.087)$} & \makecell{$4.05$ \\ $(0.60, 7.28)$} & \makecell{$0.29$ \\ $(0.09, 0.46)$} & \makecell{$\mathbf{0.013}$ \\ $(\mathbf{0.002}, \mathbf{0.031})$} \\ 
RTO-TC &\makecell{$0.857$ \\ $(0.222, 1.922)$} & \makecell{$4.72$ \\ $(0.86, 12.46)$} & \makecell{$0.30$ \\ $(0.09, 0.56)$} & \makecell{$\mathbf{0.015}$ \\ $(\mathbf{0.001}, \mathbf{0.043})$}& -- -- & -- -- & -- -- & -- -- \\ 
RTO-ODE &\makecell{$1.200$ \\ $(0.490, 2.060)$} & \makecell{$6.24$ \\ $(2.57, 9.53)$} & \makecell{$0.38$ \\ $(0.19, 0.56)$} & \makecell{$0.074$ \\ $(0.046, 0.107)$}&\makecell{$1.061$ \\ $(0.355, 2.230)$} & \makecell{$5.10$ \\ $(1.91, 7.99)$} & \makecell{$0.32$ \\ $(0.14, 0.50)$} & \makecell{$0.067$ \\ $(0.037, 0.096)$} \\ 
Reference &\makecell{$0.068$ \\ $(0.036, 0.114)$} & \makecell{$0.22$ \\ $(0.07, 0.39)$} & \makecell{$0.02$ \\ $(0.01, 0.03)$} & \makecell{$0.001$ \\ $(0.001, 0.001)$}&\makecell{$0.068$ \\ $(0.036, 0.114)$} & \makecell{$0.22$ \\ $(0.07, 0.39)$} & \makecell{$0.02$ \\ $(0.01, 0.03)$} & \makecell{$0.001$ \\ $(0.001, 0.001)$} \\ 

    \end{tabular}
    \vspace{0.15in}
    \caption{Stylized inpainting study, high noise regime: Average errors in estimation of the mean and pointwise variance of the posterior distribution, and average central moment discrepancy (CMD) and maximum mean discrepancy (MMD) errors, for different BIPSDA methods. Interdecile ranges are also displayed in parentheses. A reference that shows the average errors between two sets of i.i.d. samples from the ground truth posterior is also included. As can be seen, the `Lang' variants lead to poor estimates of the mean and variance of the posterior, while the `MAP' and `RTO' variants provide relatively strong performance. }
    \label{tab:inpaintinghighn}
\end{table*}

\section{RESULTS}
\label{sec:results}

In this section, representative results from the stylized inpainting, x-ray, and phase retrieval studies are shown to illustrate the performance of the proposed framework. 

\begin{figure*}
\includegraphics[width=\textwidth]{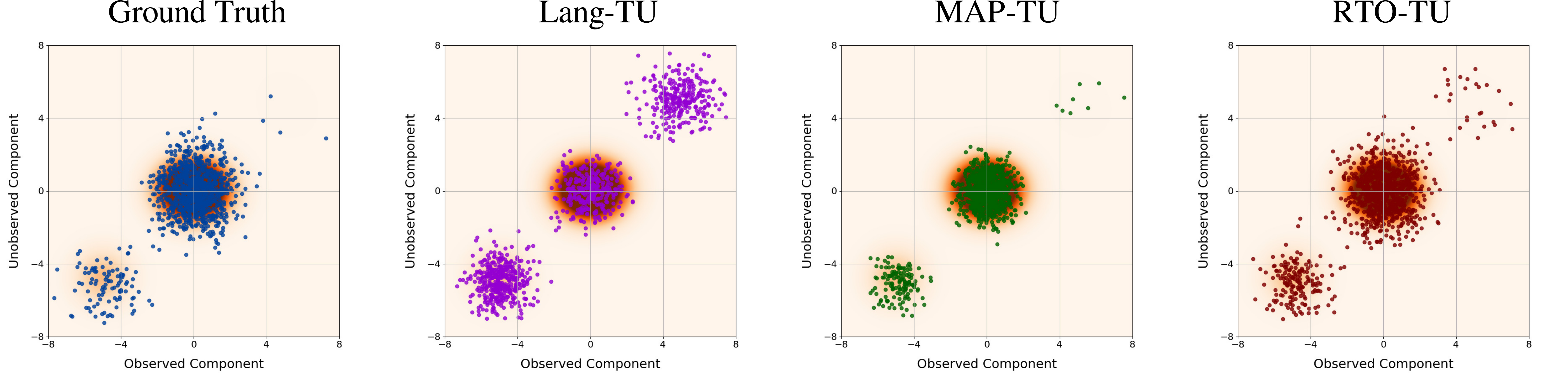}
\caption{Stylized inpainting study, high noise regime: Scatter plots of samples from three different BIPSDA methods, as well as ground truth posterior samples, in a representative trial. Two components of the ten-dimensional samples are visualized, with one of the two components in the null space of the forward operator. The samples are overlaid over the posterior density. As can be seen, the `RTO' and `MAP' variants both provide competitive performance, whereas the `Lang' variant significantly overestimates the posterior variance. }
\label{fig:inpainting_scatter}
\end{figure*}

\subsection{STYLIZED INPAINTING STUDIES}

We first examine the performance of the framework in the stylized inpainting studies. Table \ref{tab:inpaintinglown} shows the average error in the low noise regime with regards to the mean and pointwise variance of the posterior, as well as CMD and MMD metrics, for the nine BIPSDA algorithms tested. The corresponding results in the high noise regime are given in Table \ref{tab:inpaintinghighn}. Here we first note that when the analytic scores are used, the Tweedie-Correlated (`TC') based approaches perform comparably with the Tweedie-Uncorrelated (`TU') based approaches. However, `TC' variants were not implemented when using the learned scores, as the Jacobian of the learned score is not a reliable estimator of the denoising distribution covariance matrix. The performance of the `TC' variants in the learned score regime was thus omitted from this and all subsequent studies. 

All of the methods that were tested perform fairly well in the low noise regime using the learned score, with the `RTO-TU' approach providing particularly strong performance. 
In contrast, in the high noise regime, the Langevin dynamics-based approaches all perform significantly worse than the other approaches, while the `MAP' variants perform well despite lacking the theoretical basis of the `Lang' and `RTO' variants. Further, this discrepancy is not due to score modeling error, as the `Lang' variants perform poorly even when the analytic scores are used. 

Figure \ref{fig:inpainting_scatter}, which displays scatter plots of samples from three different BIPSDA methods in a representative high-noise trial, provides insight into the failure of the `Lang' variants. In particular, it can be seen that the Langevin dynamics based approach produces samples in low-density regions with respect to the likelihood and therefore overestimates the variance of the posterior distribution. The `RTO' variants, while still overestimating the posterior variance, are much more accurate. Finally, the `MAP' variants, while underestimating the variance of each of the modes, provide the most accurate estimate of the weights of the modes of all the methods tested. 

\subsection{STYLIZED X-RAY TOMOGRAPHY STUDY}

We now analyze the results from the stylized x-ray tomography study. Table \ref{tab:xray} reports the average mean, variance, CMD, and MMD errors for the different BIPSDA algorithms in this problem setting. As can be seen, in this problem setting all of the approaches we considered performed fairly well, despite the non-linearity of the forward model. The `RTO-TU' approach provides particularly strong performance and is the top performing approach with respect to all metrics we considered in the learned score regime. This demonstrates the potential of the RTO-based approaches in the context of non-linear inverse problems. 

\begin{figure*}
\includegraphics[width=\textwidth]{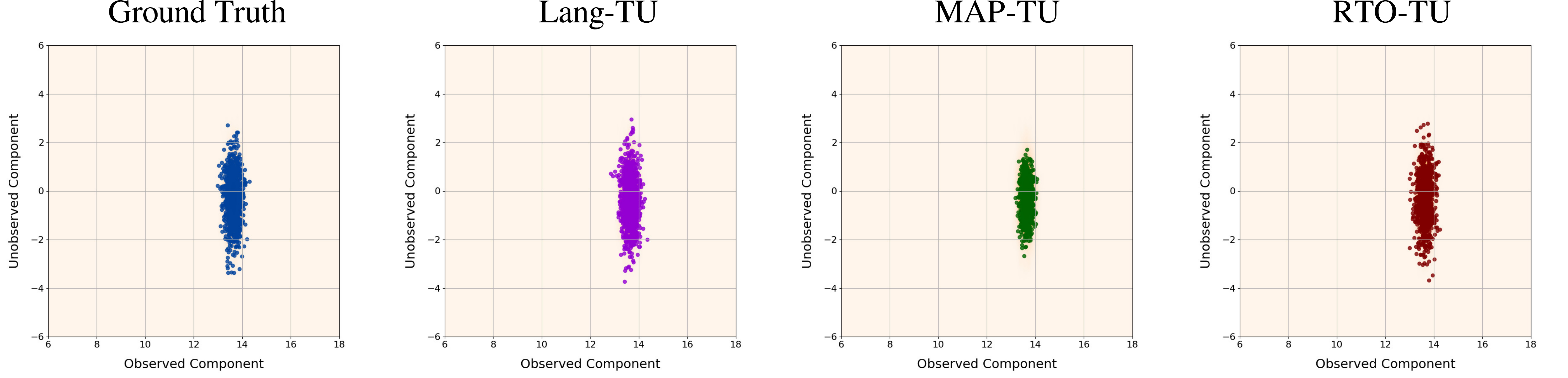}
\caption{Stylized x-ray tomography study: Scatter plots of samples from three different BIPSDA methods, as well as ground truth posterior samples, in a representative trial. Here the samples have been projected onto the space spanned by the two right singular vectors corresponding to the largest and smallest singular values of the forward operator $\mathbf{C}$. The samples are overlaid over a kernel density approximation of the pushforward of the posterior density under this projection. As can be seen, all three methods perform well in this setting, with the `Lang' and `RTO' variants providing particularly strong performance. }
\label{fig:xray_scatter}
\end{figure*}

Figure \ref{fig:xray_scatter} shows scatter plots of samples from the ground truth posterior and the `Lang-TU', `MAP-TU', and  `RTO-TU' algorithms for a representative trial. As can be seen, for this problem the posterior is approximately unimodal. Further, while all of the tested algorithms perform fairly well, the `MAP-TU' approach underestimates the posterior variance, which is consistent with the results of the stylized inpainting study. The `RTO-TU' approach is able to correct the variance underestimation of the `MAP-TU' approach and performs similarly to the `Lang-TU' method on this trial. 

\begin{table*}
\centering
    \begin{tabular}{c|cccc| cccc}
    \rowcolor{my_gray} &  \multicolumn{4}{c |}{Analytic Score} &  \multicolumn{4}{c}{Learned Score}  \\ 
    \rowcolor{my_gray} & \colorbox{my_gray}{\makecell{Mean \\ Error}}  & \colorbox{my_gray}{\makecell{ Variance \\ Error}}  &  CMD   & MMD & \colorbox{my_gray}{\makecell{Mean \\ Error}}  & \colorbox{my_gray}{\makecell{ Variance \\ Error}}  &  CMD   & MMD  \\ \hline

Lang-TU &\makecell{$0.44$ \\ $(0.27, 0.64)$} & \makecell{$0.22$ \\ $(0.04, 0.58)$} & \makecell{$0.12$ \\ $(0.07, 0.18)$} & \makecell{$0.016$ \\ $(0.008, 0.025)$}&\makecell{$0.44$ \\ $(0.26, 0.68)$} & \makecell{$0.25$ \\ $(0.05, 0.61)$} & \makecell{$0.12$ \\ $(0.07, 0.19)$} & \makecell{$0.017$ \\ $(0.009, 0.026)$} \\ 
Lang-TC &\makecell{$0.45$ \\ $(0.28, 0.65)$} & \makecell{$0.23$ \\ $(0.04, 0.60)$} & \makecell{$0.12$ \\ $(0.08, 0.18)$} & \makecell{$0.017$ \\ $(0.008, 0.028)$}& -- -- & -- -- & -- -- & -- -- \\ 
Lang-ODE &\makecell{$0.26$ \\ $(0.10, 0.48)$} & \makecell{$1.19$ \\ $(0.82, 1.57)$} & \makecell{$0.15$ \\ $(0.08, 0.23)$} & \makecell{$0.031$ \\ $(0.026, 0.036)$}&\makecell{$0.27$ \\ $(0.10, 0.51)$} & \makecell{$1.14$ \\ $(0.82, 1.43)$} & \makecell{$0.15$ \\ $(0.08, 0.22)$} & \makecell{$0.031$ \\ $(0.024, 0.037)$} \\ 
MAP-TU &\makecell{$0.12$ \\ $(0.06, 0.25)$} & \makecell{$1.66$ \\ $(1.07, 2.20)$} & \makecell{$0.14$ \\ $(0.08, 0.21)$} & \makecell{$0.162$ \\ $(0.151, 0.174)$}&\makecell{$\mathbf{0.13}$ \\ $(\mathbf{0.06}, \mathbf{0.25})$} & \makecell{$1.67$ \\ $(1.06, 2.24)$} & \makecell{$0.14$ \\ $(0.08, 0.21)$} & \makecell{$0.163$ \\ $(0.149, 0.175)$} \\ 
MAP-TC &\makecell{$\mathbf{0.09}$ \\ $(\mathbf{0.03}, \mathbf{0.26})$} & \makecell{$1.65$ \\ $(1.06, 2.20)$} & \makecell{$0.13$ \\ $(0.07, 0.22)$} & \makecell{$0.160$ \\ $(0.149, 0.171)$}& -- -- & -- -- & -- -- & -- -- \\ 
MAP-ODE &\makecell{$0.55$ \\ $(0.35, 0.75)$} & \makecell{$0.98$ \\ $(0.67, 1.26)$} & \makecell{$0.19$ \\ $(0.14, 0.25)$} & \makecell{$0.075$ \\ $(0.043, 0.112)$}&\makecell{$0.54$ \\ $(0.35, 0.70)$} & \makecell{$1.01$ \\ $(0.67, 1.33)$} & \makecell{$0.19$ \\ $(0.15, 0.25)$} & \makecell{$0.076$ \\ $(0.045, 0.113)$} \\ 
RTO-TU &\makecell{$0.12$ \\ $(0.05, 0.27)$} & \makecell{$0.19$ \\ $(0.04, 0.58)$} & \makecell{$\mathbf{0.04}$ \\ $(\mathbf{0.02}, \mathbf{0.11})$} & \makecell{$0.003$ \\ $(0.001, 0.007)$}&\makecell{$\mathbf{0.13}$ \\ $(\mathbf{0.05}, \mathbf{0.26})$} & \makecell{$\mathbf{0.19}$ \\ $(\mathbf{0.04}, \mathbf{0.55})$} & \makecell{$ \mathbf{0.04}$ \\ $(\mathbf{0.02}, \mathbf{0.10})$} & \makecell{$\mathbf{0.003}$ \\ $(\mathbf{0.001}, \mathbf{0.006})$} \\ 
RTO-TC & \makecell{$\mathbf{0.09}$ \\ $(\mathbf{0.03}, \mathbf{0.27})$} & \makecell{$\mathbf{0.18}$ \\ $(\mathbf{0.04}, \mathbf{0.58})$} & \makecell{$\mathbf{0.04}$ \\ $(\mathbf{0.01}, \mathbf{0.11})$} & \makecell{$\mathbf{0.002}$ \\ $(\mathbf{0.001}, \mathbf{0.007})$} & -- -- & -- -- & -- -- & -- -- \\ 
RTO-ODE &\makecell{$0.55$ \\ $(0.36, 0.73)$} & \makecell{$1.51$ \\ $(0.85, 2.25)$} & \makecell{$0.24$ \\ $(0.16, 0.33)$} & \makecell{$0.054$ \\ $(0.043, 0.068)$} &\makecell{$0.52$ \\ $(0.34, 0.69)$} & \makecell{$1.41$ \\ $(0.84, 2.00)$} & \makecell{$0.23$ \\ $(0.16, 0.30)$} & \makecell{$0.051$ \\ $(0.038, 0.069)$} \\ 
Reference &\makecell{$0.05$ \\ $(0.04, 0.07)$} & \makecell{$0.07$ \\ $(0.04, 0.11)$} & \makecell{$0.02$ \\ $(0.01, 0.03)$} & \makecell{$0.001$ \\ $(0.001, 0.001)$}&\makecell{$0.05$ \\ $(0.04, 0.07)$} & \makecell{$0.07$ \\ $(0.04, 0.11)$} & \makecell{$0.02$ \\ $(0.01, 0.03)$} & \makecell{$0.001$ \\ $(0.001, 0.001)$} \\ 
    \end{tabular}
    \vspace{0.15in}
    \caption{Stylized x-ray tomography study: Average errors in estimation of the mean and pointwise variance of the posterior distribution, and average central moment discrepancy (CMD) and maximum mean discrepancy (MMD) errors, for different BIPSDA methods. Interdecile ranges are also displayed in parentheses. A reference that shows the average errors between two sets of i.i.d. samples from the ground truth posterior is also included. While all methods provide competitive performance, the `RTO-TU' approach outperforms all other approaches in almost every metric. }
    \label{tab:xray}
\end{table*}

\subsection{STYLIZED PHASE RETRIEVAL STUDY}

\begin{table*}
\centering
    \begin{tabular}{c|cccc| cccc}
    \rowcolor{my_gray} &  \multicolumn{4}{c |}{Analytic Score} &  \multicolumn{4}{c}{Learned Score}  \\ 
    \rowcolor{my_gray} & \colorbox{my_gray}{\makecell{Mean \\ Error}}  & \colorbox{my_gray}{\makecell{ Variance \\ Error}}  &  CMD   & MMD & \colorbox{my_gray}{\makecell{Mean \\ Error}}  & \colorbox{my_gray}{\makecell{ Variance \\ Error}}  &  CMD   & MMD  \\ \hline

Lang-TU &\makecell{$2.42$ \\ $(0.08, 5.88)$} & \makecell{$6.26$ \\ $(1.62, 14.25)$} & \makecell{$0.20$ \\ $(0.01, 0.46)$} & \makecell{$0.076$ \\ $(0.002, 0.152)$}&\makecell{$2.45$ \\ $(0.07, 5.93)$} & \makecell{$5.73$ \\ $(0.67, 14.61)$} & \makecell{$0.20$ \\ $(0.01, 0.46)$} & \makecell{$0.073$ \\ $(0.001, 0.157)$} \\ 
Lang-TC &\makecell{$\mathbf{1.90}$ \\ $(\mathbf{0.05}, \mathbf{4.34})$} & \makecell{$\mathbf{4.01}$ \\ $(\mathbf{0.13}, \mathbf{9.69})$} & \makecell{$\mathbf{0.15}$ \\ $(\mathbf{0.00}, \mathbf{0.34})$} & \makecell{$0.048$ \\ $(0.001, 0.064)$}& -- -- & -- -- & -- -- & -- -- \\ 
Lang-ODE &\makecell{$2.53$ \\ $(0.11, 5.95)$} & \makecell{$8.33$ \\ $(3.60, 15.93)$} & \makecell{$0.21$ \\ $(0.02, 0.47)$} & \makecell{$0.134$ \\ $(0.049, 0.230)$}&\makecell{$2.54$ \\ $(0.08, 5.92)$} & \makecell{$7.85$ \\ $(2.47, 16.20)$} & \makecell{$0.21$ \\ $(0.01, 0.47)$} & \makecell{$0.128$ \\ $(0.044, 0.217)$} \\ 
MAP-TU &\makecell{$1.97$ \\ $(0.06, 4.37)$} & \makecell{$5.88$ \\ $(1.41, 11.63)$} & \makecell{$0.16$ \\ $(0.01, 0.35)$} & \makecell{$0.129$ \\ $(0.064, 0.177)$}&\makecell{$\mathbf{1.97}$ \\ $(\mathbf{0.08}, \mathbf{4.42})$} & \makecell{$\mathbf{5.63}$ \\ $(\mathbf{1.32}, \mathbf{11.46})$} & \makecell{$\mathbf{0.16}$ \\ $(\mathbf{0.01}, \mathbf{0.36})$} & \makecell{$0.127$ \\ $(0.064, 0.169)$} \\ 
MAP-TC &\makecell{$1.99$ \\ $(0.04, 4.74)$} & \makecell{$4.82$ \\ $(1.41, 10.29)$} & \makecell{$0.16$ \\ $(0.01, 0.38)$} & \makecell{$0.129$ \\ $(0.072, 0.150)$}& -- -- & -- -- & -- -- & -- -- \\ 
MAP-ODE &\makecell{$2.05$ \\ $(0.07, 4.72)$} & \makecell{$6.55$ \\ $(2.09, 13.57)$} & \makecell{$0.17$ \\ $(0.02, 0.38)$} & \makecell{$0.069$ \\ $(0.015, 0.101)$}&\makecell{$2.06$ \\ $(0.08, 4.75)$} & \makecell{$6.24$ \\ $(1.54, 12.68)$} & \makecell{$0.17$ \\ $(0.01, 0.38)$} & \makecell{$0.070$ \\ $(0.015, 0.103)$} \\ 
RTO-TU &\makecell{$2.53$ \\ $(0.11, 6.12)$} & \makecell{$6.83$ \\ $(2.44, 13.11)$} & \makecell{$0.21$ \\ $(0.02, 0.47)$} & \makecell{$0.069$ \\ $(0.006, 0.165)$}&\makecell{$2.53$ \\ $(0.09, 6.11)$} & \makecell{$5.98$ \\ $(1.59, 12.47)$} & \makecell{$0.20$ \\ $(0.01, 0.47)$} & \makecell{$\mathbf{0.065}$ \\ $(\mathbf{0.003}, \mathbf{0.156})$} \\ 
RTO-TC &\makecell{$2.04$ \\ $(0.05, 4.85)$} & \makecell{$4.68$ \\ $(0.20, 11.29)$} & \makecell{$0.16$ \\ $(0.00, 0.37)$} & \makecell{$\mathbf{0.047}$ \\ $(\mathbf{0.001}, \mathbf{0.090})$}& -- -- & -- -- & -- -- & -- -- \\ 
RTO-ODE &\makecell{$2.57$ \\ $(0.16, 5.94)$} & \makecell{$8.64$ \\ $(4.03, 14.63)$} & \makecell{$0.22$ \\ $(0.03, 0.47)$} & \makecell{$0.124$ \\ $(0.055, 0.210)$}&\makecell{$2.58$ \\ $(0.12, 5.99)$} & \makecell{$7.71$ \\ $(3.15, 13.69)$} & \makecell{$0.21$ \\ $(0.02, 0.47)$} & \makecell{$0.117$ \\ $(0.048, 0.207)$} \\ 
Reference &\makecell{$0.24$ \\ $(0.05, 0.50)$} & \makecell{$0.79$ \\ $(0.08, 1.15)$} & \makecell{$0.02$ \\ $(0.00, 0.04)$} & \makecell{$0.002$ \\ $(0.001, 0.002)$}&\makecell{$0.24$ \\ $(0.05, 0.50)$} & \makecell{$0.79$ \\ $(0.08, 1.15)$} & \makecell{$0.02$ \\ $(0.00, 0.04)$} & \makecell{$0.002$ \\ $(0.001, 0.002)$} \\

    \end{tabular}
    \vspace{0.15in}
    \caption{Stylized phase retrieval study: Average errors in estimation of the mean and pointwise variance of the posterior distribution, and average central moment discrepancy (CMD) and maximum mean discrepancy (MMD) errors, for different BIPSDA methods.  Interdecile ranges are also displayed in parentheses, and a reference that shows the average errors between two sets of i.i.d. samples from the ground truth posterior is included. Note that for the `Lang-TU' and `RTO-TU' methods in the learned score regime, 29 and 185 samples, respectively, of the million total samples generated across the hundred trials numerically diverged and were discarded. As can be seen, all of the algorithms we analyzed produce significant errors in posterior sampling on this challenging problem.}
    \label{tab:phase}
\end{table*}

\begin{figure*}
\includegraphics[width=\textwidth]{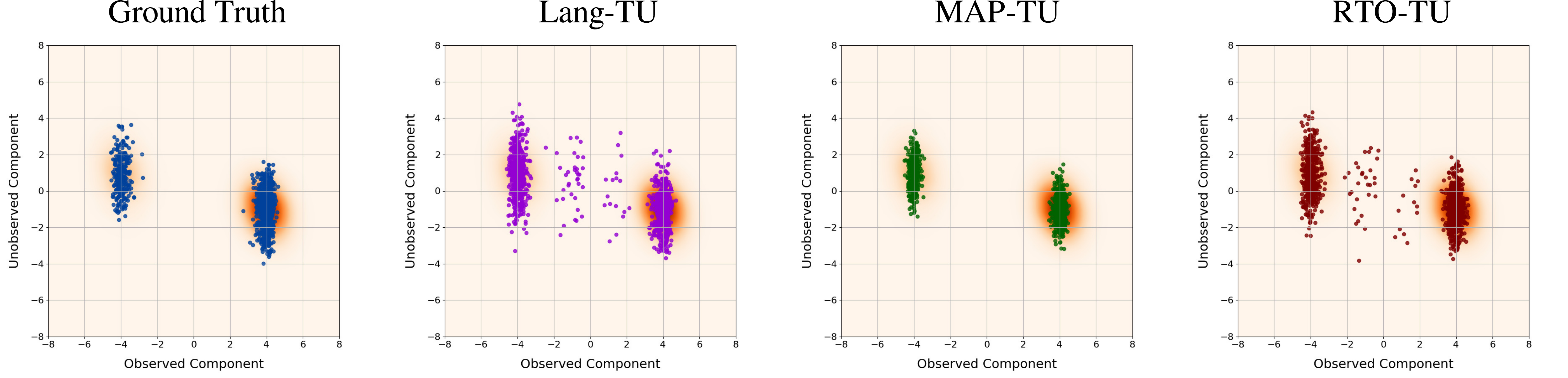}
\caption{Stylized phase retrieval study: Scatter plots of samples from three different BIPSDA methods, as well as ground truth posterior samples, in a representative trial. Here the samples have been projected onto onto the space spanned by the two right singular vectors corresponding to the largest and smallest singular values of the forward operator $\mathbf{B}$. The samples are overlaid over a kernel density approximation of the pushforward of the posterior density under this projection. As can be seen, the `RTO' and `Lang' variants both abandon samples between the two distribution modes, while the `MAP' variant underestimates the variance of the modes. }
\label{fig:phase_scatter}
\end{figure*}

We now examine the performance of the algorithms in the stylized phase retrieval study. Figure \ref{fig:phase_scatter} shows scatter plots of samples from the ground-truth posterior and the `Lang-TU', `MAP-TU', and `RTO-TU' algorithms in a representative trial (with learned score). As can be seen, both the `Lang-TU' and `RTO-TU' algorithms produce samples that lie in a low-density region between the two posterior modes. The `MAP-TU' approach avoids this issue, but suffers from significant underestimation of the variance of each of the modes.

Table \ref{tab:phase} shows the average errors in estimation of the mean and variance, as well as average CMD and MMD errors, for all nine BIPSDA methods tested. Here we first note that, when the score is known analytically, the `TC' variants perform better than the `TU' and `ODE' variants, and on some trials produce strong performance, as evidenced by the provided interdecile ranges. This indicates that the `TC' variants have some promise if a reliable approximation of the second-order score can be estimated using, e.g. the methods in \cite{meng2021estimating}. However, there is wide variance in performance across the trials, and in general the `TC' variants cannot fully overcome the issues observed with the `TU' variants in Figure \ref{fig:inpainting_scatter}. All of the BIPSDA algorithms we tested thus face performance issues in this problem context. 

\subsection{COMPUTATIONAL CONSIDERATIONS}

Here we examine the computational cost associated with the different BIPSDA algorithms. Table \ref{tab:compcost} shows the runtimes associated with each tested BIPSDA algorithm across all four studies (learned score regime). For each test-case, we measured the runtime to generate 10,000 samples per trial. The reported runtimes are the average over 100 trials. As can be seen, `TU' variants are consistently faster than their `ODE' counterparts due to the reduced number of score evaluations associated with `TU' variants. However, the dominant computational cost in our problem settings is sampling from the approximate prediction distribution, not approximating the denoising distribution. Here we find that the `MAP' and `RTO' variants are consistently faster than the `Lang' variants. The difference in runtimes is particularly pronounced in the inpainting studies, where the MAP point can be computed efficiently in closed form. Note that the cost of the Langevin dynamics based approaches could be reduced by reducing the number of Langevin iterations, but this may hurt the quality of the generated samples. In general, the computational feasibility of BIPSDA algorithms in a given problem setting is dependent on finding efficient MAP solvers (for the `MAP' and `RTO' variants) or efficient MCMC algorithms (for the `Lang' variants) for the given problem.

\begin{table*}
\centering
    \begin{tabular}{c| cccc}
    \rowcolor{my_gray} & Inpainting (Low Noise) & Inpainting (High Noise) & \; X-ray Tomography \; & Phase Retrieval \; \\

Lang-TU & $26.698 \; (0.879)$  & $26.938 \; ( 1.105$) & $514.829 \; (11.383)$ &  $333.764 \; (20.521) $\\ 
Lang-ODE & $31.698 \; (0.457)$ & $31.696 \; (0.259)$ & $519.479 \; (10.768)$ & $337.544 \; (19.719)$ \\ 
MAP-TU & $\mathbf{1.217 \; (0.011)}$ & $\mathbf{1.213 \; (0.011)}$ & $\mathbf{161.152 \; (29.001)}$ & $\mathbf{201.837 \; (32.167)}$ \\ 
MAP-ODE & $5.590 \; (0.031)$ &  $5.569 \; (0.033)$ & $\mathbf{174.486 \; (56.858)}$ & $\mathbf{210.563 \; (35.067)}$ \\ 
RTO-TU & $1.381 \; (0.002)$ & $1.409 \; (0.008)$ & $\mathbf{171.669 \; (23.894)}$ & $269.171 \; (26.925)$ \\ 
RTO-ODE & $5.780 \; (0.008)$ & $5.763 \; (0.027)$ & $\mathbf{182.396 \; (58.243)}$ & $273.492 \; (24.678)$ \\ 

\end{tabular}
\vspace{0.15in}
\caption{Computational time (in seconds) to generate $10,000$ posterior samples (averaged over $100$ trials) with different BIPSDA methods across all four studies (learned score regime). The standard deviation of the runtime is also displayed in parenthesis. All experiments were run on a Nvidia A100 GPU with 80 GB of memory. As can be seen, while computational cost is highly dependent on the specific implementation of each algorithm, in general the `MAP' and `RTO' variants are faster than the `Lang' variants. The speedup is particularly pronounced in the inpainting studies, where the MAP problem can be solved efficiently in closed form.}
\label{tab:compcost}
\end{table*}

\section{DISCUSSION AND CONCLUSION}
\label{sec:discussion_conc}

Diffusion annealing based inverse problem solvers are an emerging class of algorithms for solving Bayesian inverse problems using a diffusion model trained on the prior distribution. In this work, we introduced a general framework, Bayesian Inverse Problem Solvers through Diffusion Annealing (BIPSDA), that provides an unified formulation for this class of algorithms. This framework has two key design choices---the denoising distribution approximation and the sampler for the prediction distribution---that can be ``mixed and matched'' in a flexible manner. The large algorithmic design space of BIPSDA includes the previously introduced DAPS \cite{zhang2024improving} and DiffPIR \cite{zhu2023denoising} algorithms, which have both achieved state-of-the-art performance on different image reconstruction problems. Novel approaches can be unveiled through combinations of ideas from DAPS and DiffPIR, as well as the novel techniques proposed here to approximate the denoising distribution and to sample from the prediction distribution. Specifically, the Tweedie Correlated (`TC') technique, previously introduced in the context of the hijacking class of diffusion-based Bayesian inverse problem solvers \cite{boys2023tweedie}, is a novel contribution in the diffusion annealing context. This approach provides a theoretically-sound, hyperparameter free Gaussian approximation to the denoising distribution by leveraging the second moments of the data distribution and the generalized Tweedie formula \cite{meng2021estimating}. The randomize-then-optimize (`RTO') technique, which was originally developed as a proposal distribution for Markov chain Monte Carlo methods \cite{bardsley2014randomize,ba2022randomized,Kainan2018randomized}, is also introduced here for the first time to generate fast, approximate samples from the prediction distribution by use of an ``off-the-shelf'' numerical optimization method.

We also proposed a suite of four benchmark problems (inspired by image inpainting, x-ray tomography, and phase retrieval) for the rigorous assessment of the proposed BIPSDA framework, as well as other diffusion-based posterior samplers. A key feature of these benchmark problems is that they not only have an analytically-known posterior distribution, but also feature an analytical close-form expression for the noisy prior score. This enables analysis of the robustness of algorithms to errors in the score model by examining the idealized scenario of a known-exactly prior score.

\subsection{RESULTS ANALYSIS}

The results provide insights into both the impact of specific design choices within our framework and the performance of diffusion-model-based inverse problem solvers generally. Regarding the choice of the denoising distribution approximation, we note that approaches that use the `TC' variant provide the best performance on the two non-linear inverse problems we tested when the prior score is assumed to be known analytically. However, in our implementation this technique is currently only applicable when the score is known exactly. The training of an auxiliary neural network model to approximate higher-order moments of the data distribution is required in the learned score setting \cite{meng2021estimating}. It is also worth noting that in general, the `TU' variant
performs better than the `ODE' variant.

Regarding the sampler for the prediction distribution, the results demonstrate that the MAP estimation based approaches, despite lacking firm theoretical foundations, perform well in recovering the global structure of the posterior distribution. However, they systematically underestimate the variance of each of the posterior modes, which is unsurprising given that the `MAP' variants do not actually sample from the prediction distribution. The `RTO' technique, which is equivalent to exact sampling from the prediction distribution when the likelihood function is linear-Gaussian \cite{bardsley2012mcmc}, partially resolves this issue and provides better performance than the `MAP' variants on the stylized x-ray tomography problem and inpainting problems. Finally, the `Lang' variants work well when the posterior is unimodal (low noise regime of the stylized inpainting problem and stylized x-ray tomography problem). However, in the case of multimodal posteriors (as in the high noise regime of the stylized inpainting problem), the `Lang' variants struggle to properly incorporate the measurement information. Further, on the phase retrieval problem, Langevin dynamics without Metropolis correction completely fails, and the incorporation of a Metropolis adjustment and preconditioning was required to obtain competitive performance.

Overall, the results demonstrate that the BIPSDA framework can provide strong performance on problems with unimodal posterior distributions. Strong performance is also attainable on problems that have multi-modal posteriors and linear forward models, although in this setting the performance of the framework is sensitive to the algorithmic design choices made. On the stylized phase retrieval problem, however, all of the BIPSDA algorithms we tested produced inaccurate uncertainty estimates. This is reflective of the extremely challenging nature of the problem, for which even conventional MCMC algorithms that require knowledge of the posterior density struggle to perform well. 

\subsection{LIMITATIONS AND FUTURE DIRECTIONS}

The proposed benchmark problems provide useful insight on the ability of diffusion annealing based inverse problem solvers to provide rigorous estimates of the uncertainty of the posterior distribution. While results are promising, there are a few limitations of the present work that require further investigation. 

First, the performance of BIPSDA algorithms, and diffusion-annealing approaches in general, strongly depends on both the choice of algorithm-specific hyperparameters and the accuracy of the underlying pretrained diffusion model. In our work, we carefully controlled for these effects by fixing some of the hyperparameters (e.g. the same sequence of annealing time steps $[t_{N_A}, t_{N_A - 1}, \cdots, t_0]$ for the outer loop of Algorithm \ref{alg:bipsda}) and by comparing performances using both learned and analytic noisy prior scores. Furthermore, some of the novel algorithms that we propose within the BIPSDA framework achieve similar or superior performance to state-of-the-art methods, such as DAPS, while drastically reducing the number of hyperparameters that the user must provide. Specifically, the `TC' variant, which our work is the first to explore in the contest of diffusion-annealing approaches, provides a hyperparameter-free approximation of the denoising distribution, while on the contrary the `ODE' variant used by DAPS is highly sensitive to the choice of the discretization parameters of the probability flow ODE. Additionally, the `RTO' variant that we proposed here for the first time computes approximate samples from the prediction distribution using ``off-the-shelf'' deterministic optimization algorithms. Conversely, the `Lang' variant used by DAPS requires careful tuning of the time step size and number of time steps in the Langevin dynamics. Nevertheless, while in this work we followed recommendations in the literature when choosing hyperparameters  \cite{zhang2024improving,dhariwal_2021_diffusion}, more work is required to understand the full impact of user provided hyperparameter choices (such as the noise annealing schedule) and how to optimize the accuracy-computational cost trade-off.

Second, the algorithms in our framework were systematically tested on stylized model problems that enable principled performance analysis but have low dimensionality.  In the Supplementary Material we also  apply algorithms from the BIPSDA family to imaging problems of practical relevance. Proof-of-principle numerical results related to image inpainting (c.f. Fig. S.5 and Table S.1 in the Supplementary Material) demonstrate the computational feasibility of applying algorithms from the BIPSDA family at scale and their ability to produce plausible, diverse samples. However, to systematically evaluate the performance BIPSDA framework in these high-dimensional settings, there remains the need for image reconstruction-relevant, large-scale benchmarks with well-characterized image priors. 

Third, it is also of interest to improve the performance of BIPSDA algorithms on challenging non-linear problems like phase retrieval. Here we first note that in this study there were cases where posterior samples numerically diverged due to large error in the learned score in low probability regions of the prior. This could potentially be addressed through the incorporation of second-order score model derivative information into the score training, which has been shown to improve model accuracy in low-density regions in other problem contexts \cite{o2024derivative}. Further, even with the score known analytically, the BIPSDA algorithms performed poorly on some of the trials, and additional algorithmic innovations are required to address this. In particular, it is of interest to explore modifications to the proposed `RTO' technique, such as apodizing the noise perturbation in early BIPSDA iterations or metropolizing the RTO sampling, to improve performance on challenging problems like this one.

Finally, while the focus of this work was to provide numerical insight on the ability of decoupled noise annealing type approaches to accurately sample from the posterior distribution, the three benchmark problems that we designed are expected to become a useful tool to the community to develop, refine, and rigorously assess existing and novel approaches (including diffusion-type posterior samplers that lie outside of the BIPSDA framework), for solving Bayesian inverse problems with data-driven priors.

\begin{appendices}

\section{RELATIONSHIP BETWEEN TWEEDIE FORMULA AND THE PROBABILITY FLOW ODE}
\label{sec:appendix}

In this appendix we analyze the relationship between the  Tweedie's formula and probability flow ODE based approaches for estimation the mean of the denoising distribution $\mathbb{E}_{0 \mid t}[\bsm(0) \mid \bsm(t)]$. 

In the Tweedie's formula based approach, the estimate of the predicted mean, denoted $\bsm_{\mathrm{tweedie}}$, is computed using the pretrained score model as 
\begin{equation}
\label{eq:tweedie_pred}
\bsm_{\mathrm{tweedie}} = \bsm(t) + \sigma^2(t)  s_{\bstheta^*}(\bsm(t), t),
\end{equation}
where we have again assumed $\sigma(0) = 0$.
Since the conditional distribution of $\bsm(t)$ given $\bsm(0)$ is Gaussian, Tweedie's formula is exact and 
\begin{equation*}
\mathbb{E}_{0 \mid t}[\bsm(0) \mid \bsm(t)] = \bsm(t) + \sigma^2(t)  \nabla_{\bsm(t)} \log \pi_t(\bsm(t)). 
\end{equation*}
The only source of error in \eqref{eq:tweedie_pred} is therefore score modeling error.

In the probability flow ODE based approach, the estimated mean is obtained by solving the probability flow ODE (Eq. \eqref{eq:prob_flow_ODE}) backwards in time. For concreteness, here we consider the performance of the approach when using the backward Euler method to solve the ODE, which is commonly used in the literature \cite{song2020score, song2020denoising, zhang2024improving}. In particular, we first consider the case with only a single reverse Euler step. In this setting, the estimate of the predicted mean, denoted $\bsm_{\mathrm{ode}}$, is given by 
\begin{equation}
\bsm_{\mathrm{ode}} = \bsm(t) + t \sigma(t) \dot{\sigma}(t)  s_{\bstheta^*}(\bsm(t), t).
\end{equation}
As can be seen, the ODE and Tweedie's formula based approaches both update $\bsm(t)$ in the direction of $s_{\bstheta^*}(\bsm(t), t)$. However, the magnitude of the update differs by a factor of $d = t \sigma(t) \dot{\sigma}(t) / \sigma^2(t)$,  which in general will not be equal to one, and so therefore $\bsm_{\mathrm{ode}} \neq \bsm_{\mathrm{tweedie}}$. However, under the particular choice of parameterization $\sigma(t) = t$, $d = 1$ and the updates coincide. 

The above analysis shows that the estimate of the mean of the denoising distribution given by the Tweedie's formula based approach is exact up to error in the score model. Under a particular choice of parameterization and discretization, the probability flow ODE based approach can be made to coincide with the Tweedie's formula based approach. However, in general the probability flow ODE based approach will introduce additional approximations in the mean estimate.

\end{appendices}

\section*{ACKNOWLEDGMENTS}

The authors would like to thank Drs. Sebastian Reich, Lars Ruthotto, and Nick Alger for inspiring discussions on diffusion modeling and randomize-then-optimize algorithms.

\bibliographystyle{ieeetr}

\bibliography{refs.bib}

\begin{IEEEbiographynophoto}{Evan Scope Crafts}
graduated summa cum laude from Emory University in 2019 with a B.S. in Applied Mathematics and a B.A. in Computer Science. He was subsequently granted a masters in Computational Science, Engineering, and Mathematics from the University of Texas at Austin's Oden Institute for Computational Engineering \& Sciences in 2022, where he is currently a Ph.D. candidate. His research interests include computational imaging, inverse problems, medical imaging, statistical inference, machine learning, and optimization.
\end{IEEEbiographynophoto}

\begin{IEEEbiographynophoto}
{Umberto Villa} received the B.S. and M.S. degrees in Mathematical Engineering from the Politecnico di Milano, Milan, Italy, in 2005 and 2007, respectively, and the Ph.D. degree in Mathematics from Emory University, Atlanta, GA, USA, in 2012. He is an Assistant Professor of Biomedical Engineering and core faculty the Oden Institute for Computational Engineering and Science, The University of Texas at Austin, Austin, TX, USA. His research interests lie in the computational and mathematical aspects of large-scale inverse problems, imaging science, and uncertainty quantification.
\end{IEEEbiographynophoto}

\end{document}